\documentclass[11pt,reqno]{amsart}
\usepackage{amsmath,amsfonts,enumerate, epsfig, epsf,  indentfirst, graphicx,euscript, amssymb,amsthm,amscd}

\usepackage[latin1]{inputenc}

\makeindex

 \theoremstyle{plain}
\newtheorem{proposition}{Proposition}
\newtheorem{lemma}{Lemma}

\newtheorem{remark}{Remark}
\theoremstyle{definition}
\newtheorem{definition}{Definition}
\newtheorem{theorem}{Theorem}

\newtheorem{example}[theorem]{Example}

 \newcommand{\ep}{\epsilon}

 \def\b{\mathbf}

\title[  Mean Curvature Lines ]{Lines of   Mean Curvature  on Surfaces Immersed  in $\mathbb R^3$ }

 \author[R. Garcia]{Ronaldo Garcia}

\author[J. Sotomayor]{Jorge Sotomayor}


 \thanks{The first author was partially supported by     FUNAPE/UFG.
Both   authors are fellows of  CNPq.
 This work was done under the project PRONEX/FINEP/MCT - Conv. 76.97.1080.00 - Teoria Qualitativa das
 Equa\c c\~oes Diferenciais Ordin\'arias and had the partial support of CNPq Grant 476886/2001-5.}

\begin{document}
 \maketitle

  \begin{abstract}
 {Associated to   oriented surfaces immersed in
$\mathbb R^3,$  here are studied  pairs of transversal  foliations
 with singularities, defined on  the  {\it Elliptic} region, where the Gaussian curvature
 $\mathcal K$, given by the product of the  principal curvatures ${k_1}, k_2$ of the immersion, is positive.
The  {\it leaves} of the foliations are the  {\it lines of}     $ {\mathcal M}$-{\it mean
curvature,}
 along which   the normal curvature of the immersion is given by  a  function
  $ {\mathcal M}={\mathcal M}({k_1},{k_2})\in [{k_1}, k_2]$, called a     $ {\mathcal M}$- {\it mean
curvature}, whose  properties  extend and unify  those of the
  { {\it   arithmetic} $ {\mathcal H}=({k_1}+k_2)/2$, the  {\it geometric}
    $\sqrt{\mathcal K}$    and {\it harmonic} ${\mathcal K}/{\mathcal H}=((1/{k_1} + 1/{k_2})/2)^{-1}$  \it  classical mean
curvatures}.

The {\it singularities} of the foliations  are
 the {\it umbilic points} and {\it parabolic curves}, where ${k_1} = k_2$ and  ${\mathcal K} = 0$, respectively.

Here are determined the  patterns of       $ {\mathcal M}$- {\it mean curvature lines}  near the
 {\it umbilic points},   {\it parabolic curves}  and
 $ {\mathcal M}$-{\it mean  curvature  cycles} (the periodic
 leaves of the foliations), which are structurally  stable under small
  perturbations of the immersion. The genericity of these patterns is also established.

These patterns   provide  the three essential local ingredients to establish
 sufficient
 conditions, likely to be also  necessary,  for    $ {\mathcal M}$-{\it Mean
 Curvature Structural Stability} of immersed surfaces. This constitutes    a natural unification
and
complement  for the results obtained  previously by the authors for the
 {\it  Arithmetic},  \cite {m}, {\it Asymptotic},  \cite {a1,  a2},
  {\it  Geometric}, \cite {g} and {\it Harmonic},  \cite { h},  {\it classical}  cases  of {\it  Mean Curvature Structural Stability}.}
\end{abstract}

\vskip .4cm

\section{Introduction}
\label{sec:1}

The study of families of curves defined  for immersed surfaces by their normal
curvature properties   has attracted the interest of generations of
mathematicians, among whom are Euler,  Monge, Dupin, Gauss,
Cayley, Darboux,  Gullstrand, Caratheodory and  Hamburger, to mention only a few. See   \cite{r, St} for references.

Also dealing with families of curves,  there is {\it``The  Qualitative
Theory of Differential Equations"}  initiated by Poincar\'e and consolidated
with the study   of the {\it Structural Stability and Genericity} of differential equations
in the plane and surfaces, made systematic from $1937$ to $1962$ due  to
 the seminal works of Andronov Pontrjagin and Peixoto; see \cite {ap, mp}. The basic
  ideas  of this Theory  were extended an  applied by  Gutierrez, Garcia and Sotomayor
   to principal curvature lines \cite{gs1} as well as to  other  differential equations
   of classical geometry: the  asymptotic lines \cite {a1, a2}, the  arithmetic, geometric
    and harmonic  mean curvature lines \cite {m, g, h}.

An overview  of the ensemble of the recent
   works cited above reveals  that  they  share a common ground. In fact, there
    is a neat  analogy in  purpose, problems  and methods of analysis.
The goal  of this paper is   to inquire  more deeply on their
common
 features,  possible mathematical discrepancies and limitations of the methods used so far.

In  principle any expression such as $\mathcal M$= $\mathcal M(k_1,k_2)\in [k_1,k_2]$, involving
 the principal curvatures, could be rightly called a {\it ``mean curvature"}.
 The solutions of the differential equation: $k_n([du,dv])= \mathcal M$,
  would be called the {\it lines of}$\mathcal M$-{\it mean curvature}.  Here,
    $k_n([du,dv])= II(du,dv)/I(du,dv)$ is the normal curvature in the
    direction $([du,dv])$, as the quotient  of the {\it second } and {\it first fundamental forms} of an immersed surface.

The situations that appear in the works quoted above correspond to  the
 Principal Curvatures: $\mathcal M$=$k_1$ or $\mathcal M$=$k_2$ as well
 as to  the    Arithmetic,   Geometric and Harmonic  Mean Curvatures:
  $\mathcal M$ = $\mathcal H = ({k_1}+k_2)/2$,
  $\mathcal M$ = $\mathcal K^{1/2}$, with ${\mathcal K}= {k_1}k_2$, and  $\mathcal M$ = $\frac{\mathcal K}{\mathcal H }$.
   The asymptotic lines correspond to $\mathcal M=0$ and are supported
   by the {\it  hyperbolic} region of the immersion, where ${\mathcal K}<0$. To
   these five  functions we will refer to as the {\it ``classical" mean curvature functions.}
 In the work of the authors  the  corresponding five quadratic differential equations and respective integral foliations
 with singularities have been  unified   in terms of their properties of
   structural stability under small perturbations of the immersed surface which supports them. A
   related  unification, in terms of the  notion of $T-$ {\it systems}, focusing  the local form
    of the equations away from singularities,  was proposed  by Ogura in 1916; see \cite {og}.

In this work we
 extend  to a general mean curvature function  $\mathcal M$, with mild regularity assumptions,  the essential  results obtained for the  ``classical" mean curvature functions: principal lines \cite{gs1},  asymptotic lines \cite {a1, a2},  arithmetic, geometric  and harmonic  mean curvature lines \cite {m, g, h}.
This generalization  also includes interesting  cases of {\it
 mean curvature functions} which seem to have been overlooked previously in Geometry.

\vskip 0.1cm

To make precise the  requirements
imposed on a   function $\mathcal M$= $\mathcal M(k_1,k_2)$ to be
called  a
  {\it  mean curvature function}, it is appropriate to write its   expression  as $\mathcal M=m(H,K)$, in terms of a function $m$  of  the  $H,\, K$ variables, which  when  replaced by the {\it elementary symmetric functions} of $(k_1,k_2)$: $H=\mathcal H(k_1,k_2)$ and $K=\mathcal K(k_1,k_2)$ give back $\mathcal M(k_1,k_2)$.

\begin{definition}\label{def:m}
 A  function $\mathcal M =m(H,K)$ is called a mean curvature function
provided the following holds:

\begin{itemize}
\item[1)] It satisfies  $(m-H)^2 \leq H^2-K$ on  the region $H^2
\geq K   \geq 0$   (mean function condition).

 \item[2)]It  is
continuous on the region    $H^2 \geq  K  \geq 0$  and analytic on
$ K > 0$ (basic regularity condition).

\item[3)]  $ m(tH,t^2K)=t m(H,K)$, $\; t\geq 0$,  (weighted
homogeneity condition).
 \noindent

\end{itemize}
\end{definition}

 \begin{remark}\label {rem:m} In terms of $(k_1,k_2)$, definition \ref{def:m} amounts to the following:
\begin{itemize}
\item[1)]  ${\mathcal M}={\mathcal M}({k_1},{k_2})\in [{k_1},{k_2}]$, is symmetric i.e. ${\mathcal M}({k_1},{k_2})={\mathcal M}({k_2},{k_1})$,
\item[2)]  it is continuous everywhere and analytic on
 $k_1k_2 > 0$   and
\item[3)] it is homogeneous i.e. ${\mathcal M}({tk_1},{tk_2})={t\mathcal M}({k_1},{k_2})$.
\end{itemize}

Notice that $ m(H,H^2)=H$, the diagonal condition, also expressed by  ${\mathcal M}({k},{k})=k$, follows directly from the {\it mean function condition}, labeled  1).  One can also   pass from  ${\mathcal M}({k_1},{k_2})$ to  $m(H,K)$ by the transformation $k_1=H-\sqrt{H^2-K}$,  $k_2=H+\sqrt{H^2-K}$
\end{remark}

For more developments on the subject of {\it Means}, carried out
   from the perspective of Arithmetic and Analysis, the reader is addressed to Borwein and Borwein \cite {bb}, Hardy
et al. \cite{hl} and Mitrinovic \cite{M}.  Definition \ref{def:m}
is adapted for our needs from \cite {bb} , Chapter 8. Additional
requirements will be made to it later  to deal with differential
geometric problems.

  Definition \ref{def:m}  includes the classical symmetric means and also their most important  generalizations, of which we review two examples below.

\begin{example}\label {ex:1}
The    {\it Holder  mean of order r}:

$$  {\mathcal H_r}({k_1},{k_2})=(\frac{{k_1}^r+ {k_2}^r}{2})^{1/r},\; r\neq 0;\;\; {\mathcal H_0}({k_1},{k_2})=\sqrt{k_1k_2}$$

\noindent  generalizes the    classical {\it   arithmetic, geometric and  harmonic  mean
curvatures}, given respectively  by  $r=1, \, 0, \, -1$. The continuity at $r=0$ and the limits $\lim_{r\to{\pm\infty}}{\mathcal H}_r=H\pm\sqrt{H^2-K}$ are well known; see \cite {hl, M}.
\end{example}

Taking $r$ as a parameter,  this  defines a natural transition between
 the classical means and consequently between their associated differential equations and foliations with singularities.

\begin{example}\label {ex:2}

The    classical  {\it AG  mean} of Gauss and Legendre is
  defined  by

$$ {\mathcal A}{\mathcal G} ({k_1},{k_2})= I(1,1)/I({k_1},{k_2});\;I({k_1},{k_2})=\int_0^\infty\frac{dt}{(t^2+{k_1}^2)^{1/2}(t^2+{k_2}^2)^{1/2}}. $$
\end{example}

 In  Borwein and Borwein \cite {bb} can be found an enlightening study of this mean as well as a general treatment  of the basic properties of  Mean Functions.
See also Weisstein \cite {w}  for other  references  on {\it means},  including  recent, non symmetric, generalizations of the  {\it AG  mean}.

\vskip .3cm

For any  mean curvature  function $\mathcal M$,  as in definition
\ref{def:m},  are defined  two transversal foliations whose
leaves, called the {\it lines of} $\mathcal M$-{\it mean
curvature},  are  the solutions of the quadratic  differential
equation  $k_n([du,dv])= \mathcal M$. These foliations, called
here the $\mathcal M$-{\it    mean curvature foliations}, are well
defined and regular only on the non-umbilic part of the elliptic
region of the immersion, where the  Gaussian Curvature is
positive. The set where the Gaussian Curvature vanishes,
 the parabolic set, is generically a  regular curve which is the border of the elliptic region.
The umbilic points
 are those at which the principal curvatures coincide,  generically are isolated and disjoint
 from the parabolic curve.

The transversal foliations, are
assembled with the umbilic and parabolic points to define the $\mathcal M$-{\it    mean curvature configuration}
 of an immersed surface. See section 2 for precise definitions.

This paper establishes  sufficient conditions, likely to be also necessary, for
the  structural stability of $\mathcal M$-{\it    mean curvature configurations},
  under small perturbations of the immersion. See sections 2 and 6  for precise statements.

Three  local ingredients  are essential  to express these sufficient conditions:
the {\it umbilic
points}, endowed  with  their  $\mathcal M$-{\it mean curvature  separatrix structure}, the
$\mathcal M$-{\it mean curvature cycles}, with the calculation of the derivative of the Poincar\'e return map,
through which is expressed the hyperbolicity condition and   the {\it parabolic curve},
 together with the {\it parabolic tangential  singularities}
and associated {\it separatrix structure}.

The conclusions  of this work, on the {\it elliptic  region}, are
complementary to  results valid independently for {\it asymptotic
foliations} on the {\it hyperbolic region}  (on which the Gaussian
curvature is negative), for which  the {\it separatrix  structure}
near the parabolic curve and the asymptotic structural stability
has been studied in \cite {a1, a2}.

\vskip 0.2cm

This paper is organized as follows:

Section \ref{sec:2} is devoted to the general study of the differential equations and general properties of  $\mathcal M$-Mean Curvature Lines. Here are given the precise definitions of the  $\mathcal M$-Mean Curvature Configuration  and of the  two transversal  $\mathcal M$-Mean Curvature Foliations with singularities into which it splits.  The definition of  $\mathcal M$-Mean Curvature Structural Stability focusing on the preservation  of the qualitative properties of the foliations and the configuration under small perturbations of the immersion, will be given at the end of this section.

In Section \ref{sec:3} the equation of lines of  $\mathcal M$-mean curvature is written in a Monge chart. The condition  for the   $\mathcal M$-mean curvature structural stability  at umbilic points  is explicitly stated in terms of the coefficients of the third order jet of the    function which represents the immersion in a Monge chart. The local  $\mathcal M$-mean curvature separatrix configurations at stable umbilics is established for $C^4$ immersions. The patterns  resemble those established for  the three  Darbouxian umbilic points in  the stable  arithmetic mean curvature configurations \cite {da, gs1}.

  In Theorem \ref {th:31} it is proved that this is due to the  properties of mean curvature functions (definition \ref {def:m} and remark \ref {rem:m} ). This clarifies    why they appear also  in the  geometric and harmonic mean curvature configurations studied previously in \cite {g, h}.

In Section \ref{sec:4} is calculated the derivative of first return Poincar\'e map along
a  $\mathcal M$-mean curvature cycle. It consists of an integral
expression
involving  $\mathcal M$ and  other natural curvature  functions
along the  cycle. Under an additional regularity condition on $\mathcal M$ (or $m$), denominated {\it positive regularity} in definition \ref {def:bar},  it is shown how to deform an immersion so that a non hyperbolic $\mathcal M$-mean curvature cycle becomes hyperbolic.

In Section \ref{sec:5} are studied the foliations by  lines of  $\mathcal M$-
mean curvature  near the parabolic set   of an immersion, which
typically is   a regular curve.
Here it is also necessary to impose  additional regularity  conditions on the function  $\mathcal M$. Two cases are considered in detail, denominated {\it 1-regular} and {\it 1/2-regular}. See definitions \ref  {def:mb1} and \ref {def:mb2}.
In  the 1-regular (resp.  1/2-regular) case,
three (resp. only two) singular tangential patterns
exist generically:  the {\it folded node},  the  {\it
folded saddle} and the {\it folded focus} (resp. only the  {\it folded node}  and the  {\it
folded saddle}).
The results of this section  extend  those obtained for  the harmonic, as in \cite {h}, (resp. geometric, as in \cite {g},)  mean curvature.

In Section \ref{sec:6} the results presented in Sections \ref{sec:3}, \ref{sec:4} and \ref{sec:5} are put together
to provide sufficient conditions for  $\mathcal M$-Mean Curvature Structural Stability.

  The density of these conditions  is established in  section \ref{sec:7}.
The delicate point here is the elimination of non-trivial recurrent $\mathcal M$-mean curvature lines by means of small perturbations of the immersions.  The main steps for the somewhat technical proof  of this part are explained in detail here under suitable hypotheses.

Section \ref{sec:8} presents  a short overview of the
achievements of this paper and points out to  some possible lines
for  future research.

 For a discussion on  historic  grounds of the prominence of the classical means
  in Arithmetic, Geometry and Analysis and of  the needs for their generalization,
  the reader is addressed to  the
  book by Borwein and Borwein, \cite {bb}. See also the essay by Wassell, \cite {wa}.

 \section{Differential Equations for   $\mathcal M$- Mean Curvature Lines}\label{sec:2}

Let $\alpha : {\mathbb V}^2\to \mathbb R^3$ be a $C^r,\;\; r\geq 4,$ immersion of
an oriented smooth surface ${\mathbb V}^2$ into $\mathbb R^3$. This means that $D\alpha$
is injective  at every point in ${\mathbb V}^2$.

 The space $\mathbb R^3$ is oriented by a  once for all fixed orientation
and endowed with the Euclidean inner product  $<,>$.

Let $N$  be the positive unit   vector field normal to $\alpha$.
This means  that for  any  positive chart $(u,v)$ of ${\mathbb V}^2$,  $\{\alpha_u, \alpha_v, N\}$
is  a positive frame in $\mathbb R^3$.

In such   chart $(u,v)$, the {\it first fundamental form} of an immersion
$\alpha$ is
given by:

\centerline{$I_\alpha= <D\alpha,D\alpha>= E du^2+2Fdudv+Gdv^2$,}

\noindent  with $E=<\alpha_u,\alpha_u>$, $F=<\alpha_u,\alpha_v>$, $G=<\alpha_v,\alpha_v>$

The  {\it second fundamental form} is given by:

 \centerline{$II_{\alpha}= <N , D^2\alpha> = e du^2 + 2f dudv + g dv^2$,}

\noindent  with $e=<N,\alpha_{uu}>=-<N_u,\alpha_u>$, $f=<N,\alpha_{uv}>=-<N_u,\alpha_v>$, $g=<N,\alpha_{vv}>=-<N_v,\alpha_v>$.

The normal curvature  at a point $p$ in a tangent direction $t=[du:dv]$ is given by:

 $$ k_n= k_n(p) =\frac{ II_\alpha(t,t)}{I_\alpha(t,t)}.$$

Given a mean curvature  function $\mathcal M$  as in definition \ref{def:m} and remark \ref {rem:m},   the lines of $\mathcal M$-mean curvature  of $\alpha$ are regular curves  on $\mathbb V^2$  along which the normal  curvature is    equal to $\mathcal M$. That is,  $k_n={\mathcal M}({k_1},{k_2})=m({\mathcal H},{\mathcal K})$, where $ {\mathcal K}$ = ${\mathcal K}{_\alpha}$ and $ {\mathcal H}$ = ${\mathcal H}{_\alpha}$ are  the Gaussian and Arithmetic Mean curvatures of $\alpha$. See definition \ref{def:m}.

Therefore the   pertinent differential equation for these lines is given by:

$$\frac{edu^2+2fdudv+gdv^2}{Edu^2+2Fdudv+Gdv^2} =  {\mathcal M}({k_1},{k_2})=m({\mathcal H},{\mathcal K})$$

\noindent where  $ {\mathcal H}=\frac{Eg+eG-2fF}{2(EG-F^2)}$ and $ {\mathcal K}=\frac{eg-f^2}{ EG-F^2 }$ or equivalently, according to remark \ref{rem:m}, expressing  in ${\mathcal M}({k_1},{k_2})$  the principal curvatures  in terms of  $({\mathcal H},{\mathcal K})$.

Or equivalently by

 \begin{equation}\label{eq:harm}
 (g-{\mathcal M}G ) dv^2  + 2(f-{\mathcal M}F) dudv + (e-{\mathcal M}E)du^2=0.
\end{equation}

This equation is well defined
 only on the closure of the {\it Elliptic region}, ${\mathbb E}{\mathbb V}^2_\alpha$,   of $\alpha$, where ${\mathcal K} >0$. It is bivalued and $C^{r-2},\;\; r\geq 4,$ smooth on the complement of the umbilic,  ${\mathcal U}_\alpha$, and parabolic, ${\mathcal P}_\alpha$, sets of the immersion $\alpha$.  In fact, on ${\mathcal U}_\alpha$ , where the principal curvatures  coincide, i.e where ${\mathcal H}^2-{\mathcal K}=0$, the equation vanishes identically; on  ${\mathcal P}_\alpha$, it is univalued when $k_2={\mathcal M}$ or when ${\mathcal M}=k_1=0$.

\vskip 0.3cm

 The  lines of  $\mathcal M$-mean curvature of immersions  will be  organized into  the
$\mathcal M$-{\it  mean curvature
configuration,}  as follows:

 Through every point $p\in {\mathbb E}{\mathbb V}^2{_\alpha} \setminus ({\mathcal U}_\alpha \cup {\mathcal P}_\alpha)$,  pass two
 $\mathcal M$-mean curvature
lines of $\alpha$. Under the orientability
hypothesis imposed on $\mathbb V^2$, the  $\mathcal M$-mean curvature
lines define two foliations:  $\mathbb H_{\alpha,1}$,
called the $\mathcal M$-{\it minimal   mean curvature foliation}, along which the
{\it geodesic torsion} is negative (i.e  $\tau_g=- \sqrt {(k_2-{\mathcal M})({\mathcal M}-k_1)})$
 and $\mathbb H_{\alpha,2}$, called the $\mathcal M$-{\it maximal
  mean curvature foliations}, along which the geodesic torsion is
positive
(i.e  $\tau_g = \sqrt {(k_2-{\mathcal M})({\mathcal M}-k_1)})$

By comparison with the arithmetic mean curvature directions, making angle $\pi/4$ with the minimal principal directions, the $\mathcal M$ directions  are located  between them and the principal ones, making an angle $\theta_m$ such that
$\tan\theta_m$ =$\pm\sqrt{\frac{{\mathcal M}-k_1}{k_2-{\mathcal M}}}$,
 as follows from Euler's Formula: $k_n(\theta)=k_1\cos^2(\theta)+k_2\sin^2(\theta)$, which leads directly to $\cos\theta_m$ =$\pm\sqrt{\frac{{\mathcal M}-k_1}{k_2-k_1}}$ and $\sin\theta_m$ =$\pm\sqrt{\frac{k_2-{\mathcal M}}{k_2-k_1}}$. The particular expression for the geodesic torsion given above
results from the formula  $\tau_g=(k_2-k_1)\sin{\theta} \cos{\theta}$ \cite{St}. See also lemma \ref{lem:41}   in Section \ref{sec:4} below.

\vskip 0.3cm

The quadruple
$$\mathbb H_\alpha^{\mathcal M}=\{{\mathcal P}_\alpha,
{\mathcal U}_\alpha, \mathbb H_{\alpha,1}^{\mathcal M},\mathbb H_{\alpha,2}^{\mathcal M}\}$$
\noindent is called the ${\mathcal M}$-{\it    mean curvature configuration} of $\alpha$.

It splits into two foliations with singularities:
$$\mathbb G^{\mathcal M}_{\alpha,i}=\{{\mathcal P}_\alpha,
{\mathcal U}_\alpha, \mathbb H^{\mathcal M}_{\alpha,i}\},\, i=1, 2.$$

Let $\mathbb V^2$ be also   compact. Denote by ${\mathcal I}^{r,s}({\mathbb V^2)}$ the space of $ C^r $ immersions of $\mathbb V^2$ into the Euclidean space $\mathbb R^3$, endowed with the $C^s$ topology.

 An immersion $\alpha$ is said $C^s$-$\mathcal M$-{\it local   mean curvature structurally
stable at a compact set} $C\subset \mathbb V^2$
if for any sequence of immersions $\alpha_n$ converging to $\alpha$ in ${\mathcal I}^{r,s}({\mathbb V^2)}$
there is a neighborhood $V_C$ of $C$, sequence of compact subsets $C_n$ and a sequence of homeomorphisms mapping $C$
to $C_n$, converging to the identity of $\mathbb V^2$ such that on $V_C$ it  maps umbilic and parabolic  points and  arcs of the   $\mathcal M$-mean curvature foliations
$\mathbb H_{\alpha,i}$  to those  of $\mathbb H_{\alpha_n,i}$
for $i = 1,\;2$.

An immersion $\alpha$ is said to be $C^s$-$\mathcal M$-{\it     mean curvature structurally stable} if the compact $C$ above is the closure of  ${\mathbb E}{\mathbb V}^2{_\alpha}$.

Analogously, $\alpha$ is said  to be {\it i-} $C^s$-$\mathcal M$-{\it    mean curvature structurally stable}
if only the preservation of elements of  {\it i-th, i=1,2} foliation with singularities is required.

A general study of the structural stability of  quadratic
differential equations (not necessarily derived from normal
curvature properties) has been carried out by Gu{\'\i}\~nez \cite
{gn}. See also the work of Bruce and Fidal  \cite{bf} , Bruce and
Tari \cite{bt} and   Davydov \cite {dav}
 for the analysis of umbilic points  for general quadratic
 and also implicit
 differential equations.

For a study of the topology of  foliations with  non-orientable
singularities
on two dimensional manifolds,
see the works of Rosenberg and Levitt \cite{ro,le}. In these works the
leaves  are not  defined by normal  curvature properties.

\section{  $ {\mathcal M}$-mean curvature  lines near umbilic points }\label{sec:3}

 Let  $0$ be an umbilic point of a $ C^r ,\; r\geq  4,$ immersion  $\alpha$ parametrized in a Monge chart $(x,y)$ by $\alpha(x,y)=(x,y,z(x,y))$, where

\begin{equation} \label{eq:1} z(x,y)=\frac k2(x^2+y^2)+ \frac a6 x^3+\frac b2 xy^2+\frac c6 y^3+O(4)\end{equation}

This reduced form is obtained by means of a rotation of the $x,y$-axes. See \cite {gs1, gs2}.

\begin{figure}[htbp]
 \begin{center}
 \includegraphics[angle=0, width=11cm]{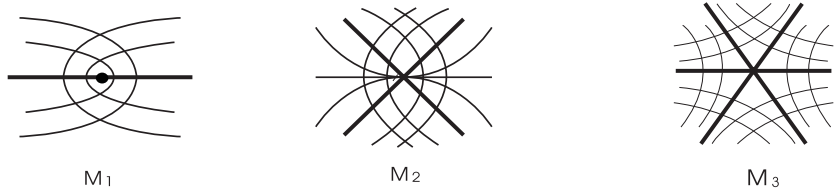}
 \caption{$ {\mathcal M}$-mean curvature lines near the umbilic points $M_i$ and their separatrices \label{fig:umb}}
   \end{center}
 \end{figure}

\begin{proposition}\label{prop:2}  Assume  the notation  established in equation \ref{eq:1}.
Suppose that the transversality condition $T_m :  kb(b-a)\ne 0$  holds
 and
consider the following situations:
\begin{itemize}
\item[$M_1$)]   $\;\;\Delta_{m}>0$

\item[$M_2$)]$\;\;\Delta_{m}<0 $ and $ \;\;  \dfrac ab  >1$

\item[$M_3$)]    $\;\; \dfrac ab < 1. $
\end{itemize}
Here  $ \Delta_{m} =4c^2(2a-b)^2-[3c^2+(a-5b)^2][3(a-5b)(a-b)+c^2].  $

  Then for every  mean curvature function $\mathcal M$, the
foliations $\{\mathbb H^{\mathcal M}_1, \mathbb H^{\mathcal M}_2\}
$ have  in a neighborhood of $0$,  one hyperbolic sector in the
$M_1$ case, one parabolic and one hyperbolic sector in $M_2$ case
and three hyperbolic sectors in the case $M_3$. These  points are
called  Darbouxian umbilics, see Figure \ref{fig:umb}.

  The {\it separatrices} of these singularities are called {\it umbilic separatrices}.
\end{proposition}

\begin{proof}
Near $0$, the functions $\mathcal K$ and $\mathcal H$ have the following Taylor expansions.

$$  {\mathcal K}=k^2+  (a+b) k x+   ck y+O_1(2),\hskip .5cm
{\mathcal H}=k+ \frac 12(a+b)x+\frac 12 cy+O_2(2).$$

Assume that
${\mathcal M}$ is a mean curvature function, as defined in \ref{def:m} and therefore    from $(m-H)^2\leq H^2-K $  follows that
 $j^1{\mathcal M}_{(H,H^2)}(0,0)=j^1{\mathcal H}_{(H,H^2)}(0,0)$.

The  differential equation of the   $ {\mathcal M}$-mean curvature lines
\begin{equation} \label{eq:mh} [g-{\mathcal M}G]dx^2+2[f-{\mathcal M}F]dxdy+  [e-{\mathcal M}E]dy^2=0 \end{equation}
\noindent  is given by:

\begin{equation}\label{eq:mg}
\aligned\; [(b-a)x+cy+R_1(x,y)]dy^2 +& [4by+R_2(x,y)]dxdy \\
\; -[(b-a)x  +  cy+R_3(x,y)]dx^2 =&0
\endaligned
\end{equation}

 \noindent where $R_i$,  $i=1,2,3$,  represent  functions of
order $O((x^2+y^2)).$

Thus, at the level of first jet, the differential equation \ref{eq:mg} is the same as that of the arithmetic mean curvature lines given by
$$ [g- {\mathcal H}G]dv^2+2[f- {\mathcal H}F]dudv+  [e- {\mathcal H}E]du^2=0.$$

The conditions on $\Delta_m $ coincide with  those on   $\Delta_H$, established to characterize the arithmetic mean curvature Darbouxian umbilics studied in detail  in \cite{m}. Thus reducing the analysis of the umbilic points to that of the  hyperbolicity of   saddles and nodes of plane vector fields, whose phase portraits are determined   only by  the first jets of the vector field at the singularities,   which are  calculated only in terms of the first jet of the equation \ref {eq:mg} at the umbilic point.
\end{proof}

\begin{theorem}\label{th:31}
An immersion $\alpha\in {\mathcal I}^{r,s}({\mathbb V^2)}$, $ r\geq   4$, is $C^3- {\mathcal M}$-local
   mean curvature structurally  stable at ${\mathcal U}_\alpha$ if and only if
every $p\in {\mathcal U}_\alpha$ is one of the types $M_i$, $i=1,2,3$ of proposition \ref{prop:2}.
\end{theorem}

\begin{proof} Clearly proposition \ref{prop:2} shows that the condition $M_i$, $i=1,2,3$ together with $T_m:kb(b-a)\ne 0 $  imply
the $C^3-$ $ {\mathcal M}$-local  mean curvature structural stability. This involves the construction of  the homeomorphism (by means of
canonical regions), mapping simultaneously  minimal and maximal   $ {\mathcal M}$-mean curvature lines  around the umbilic points of
$\alpha$ onto those of a $C^4$ slightly perturbed immersion.

We will discuss the necessity of the  condition $T_m:k(b-a)b\ne
0$ and of the conditions $M_i$, $i=1, \,2,\, 3$.  The first one  follows from
its identification with a transversality condition that guarantees the persistent isolation of the umbilic
points of $\alpha$ and its separation from the parabolic set, as well as the persistent regularity of  the Lie-Cartan
  surface
${\mathcal G}$, obtained from the projectivization of the equation \ref{eq:mh}.  Failure of
$T_m$
condition has the following implications:
\begin{itemize}
\item[a)] $b(b-a)=0$; in this case the elimination or splitting of the umbilic point can be achieved by small perturbations.

\item[b)] $k=0$ and $b(b-a)\ne 0$; in this case a small perturbation separates the umbilic point from the parabolic set.
\end{itemize}
 The necessity of condition  $M_i$ follows from its dynamic identification with the hyperbolicity  of the
equilibria along the projective line of the vector field obtained lifting equation
 (\ref{eq:mh}) to the surface   ${\mathcal G}$. Failure of this condition would
 make possible
to  change the number of  $ {\mathcal M}$-mean curvature separatrices
at the umbilic point by means a small  perturbation of
the immersion.
\end{proof}

\section{Periodic   $ {\mathcal M}$-Mean Curvature Lines}\label{sec:4}

Let $\alpha:\mathbb V^2\to \mathbb R^3$ be an immersion of a compact and
 oriented surface and consider the foliations $\mathbb H^{\mathcal M}_{\alpha, i}$, $i=1,\;2$, given by
 the $ {\mathcal M}$-{\it mean curvature lines}.

  In terms of ${\mathcal M}$ and other natural  geometric invariants of the immersion, here is established an integral
expression  for  the first derivative of the  return map of a
periodic    $ {\mathcal M}$-mean curvature line, called  also a $ {\mathcal M}$-{\it    mean curvature cycle}.
Recall that the return map associated  to a cycle
is a local diffeomorphism with a fixed point, defined on a cross section  normal to the cycle by following  the
integral curves through this section until they meet again the section. This map is called holonomy in   Foliation
Theory and Poincar\'e Map in Dynamical Systems, \cite{pm}.

A   $ {\mathcal M}$-mean curvature cycle is called {\it hyperbolic} if the first derivative of the return map at the fixed  point is
different from one.

{ The   $ {\mathcal M}$-mean curvature foliations $\mathbb H_{\alpha,i}$ have no  cycles such that  the return map reverses the orientation.}
 Initially, the integral expression for the derivative of the return
map is
obtained in class $ C^6$; see Lemma \ref{lm:42} and Proposition \ref{prop:41}.

The characterization of hyperbolicity of  $ {\mathcal M}$-mean curvature cycles in terms
of local structural stability is given in  Theorem  \ref{th:41} of
this section.

\begin{lemma}\label{lem:41} Let $c:I\to {\mathbb V}^2$ be a  $ {\mathcal M}$- mean curvature
line parametrized by arc length. Then the Darboux frame is given by:

$$\aligned
 T^\prime &= k_g N\wedge T+ {\mathcal M} N\\
(N\wedge T)^\prime &= -k_g T+ \tau_g N\\
N^\prime &= -{\mathcal M} T-\tau_g N\wedge T\endaligned\eqno
$$
 \noindent where $\tau_g=  \pm \sqrt {( {\mathcal M}-k_1)(k_2-{\mathcal M})} $.
The sign of $\tau_g$
is positive (resp. negative) if $c$ is maximal (resp. minimal)  $ {\mathcal M}$-mean curvature line.
\end{lemma}

\begin{proof} The normal curvature $k_n$ of the curve $c$ is by the definition
 the    mean curvature function ${\mathcal M}$.
>From the Euler equation $k_n=k_1\cos^2\theta+k_2\sin^2\theta={\mathcal M}$, get
$\tan\theta=\pm\sqrt{\frac{{\mathcal M}-k_1}{k_2-{\mathcal M}}}$.

  Therefore, by direct calculation, the geodesic torsion is given by
$\tau_g$=$(k_2-k_1)\sin\theta\cos\theta$  =
 $\pm \sqrt{  ({\mathcal M}-k_1)(k_2-{\mathcal M})}$.
\end{proof}

 \begin{remark}
The expression for the geodesic curvature $k_g$ will not be needed explicitly in this work. However, it can be given in terms of the principal curvatures and their derivatives using a formula due to Liouville \cite{St}, pp. 130-131.
In fact, in a principal chart $(u,v)$ the geodesic curvatures of the coordinate curves are given by:
$$ k_g|_{v=v_0}=\frac{\frac{\partial k_1}{\partial v}}{k_1-k_2}, \;\;\;
 k_g|_{u=u_0}=\frac{\frac{\partial k_2}{\partial u}}{k_1-k_2}.$$

Therefore, by Liouville formula, the geodesic curvature of a curve
$c(s)$ parametrized by arc length and  that makes an angle
$\theta(s)$ with the principal direction $e_1=\partial /\partial
u$ is
$$k_g=\frac{d\theta}{ds}+k_g|_{v=v_0}\cos\theta +k_g|_{u=u_0}\sin\theta.$$
\end{remark}

\begin{lemma}\label{lm:42} Let $\alpha:\mathbb V^2\to \mathbb R^3$ be an immersion of class $ C^r $, $ r\geq   6$, and $c$ be a $\mathcal M$ - mean
curvature cycle of $\alpha$,
 parametrized by arc length and  of length $L$. Then the expression,

$$\alpha(s,v)=c(s)+ v(N\wedge T)(s)+[ (2{\mathcal H} (s)- {\mathcal M}(s))   \frac{v^2}2+
 \frac{A(s)}{6}v^3+ v^3B(s,v)]N(s)\eqno$$

 \noindent where $B(s,0)=0$, defines    a local chart $(s,v)$ of class $C ^{r-5}$ in a neighborhood of $c$.
\end{lemma}

\begin{proof}  The curve $c$ is of class $C ^{r-1}$ and the map
$\alpha(s,v,w)=c(s)+ v(N\wedge T)(s)+wN(s)$ is of class $C ^{r-2}$ and is a local diffeomorphism  in a neighborhood of the axis $s$. In fact $[\alpha_s,\alpha_v,\alpha_w](s,0,0)=1$. Therefore there is a function $W(s,v)$ of class $C ^{r-2}$ such that $\alpha(s,v,W(s,v))$ is a parametrization of a tubular neighborhood of $\alpha\circ c$. Now for each $s$, $W(s,v)$ is just a parametrization of
the curve of intersection between $\alpha(\mathbb V^2)$ and  the normal plane generated by $\{(N\wedge T)(s), N(s)\}$. This curve of intersection  is tangent to $(N\wedge T)(s)$ at $v=0$ and notice that  $k_n(N\wedge T)(s)=2{\mathcal H}(s)-{\mathcal M}(s) $. Therefore,

\begin{equation} \aligned  \alpha(s,v,W(s,v))=& c(s)+ v(N\wedge T)(s)\\
+&[(2{\mathcal H} (s)- {\mathcal M}(s))\frac{ v^2}2+
 \frac{A(s)}6v^3+ v^3B(s,v)]N(s), \endaligned \end{equation}

 \noindent where $A$ is of class $C ^{r-5}$ and $B(s,0)=0$.
\end{proof}

We now  compute the coefficients of the first and second fundamental forms in  the chart $(s,v)$ constructed above, to be used in  proposition \ref{prop:41}.

$$\aligned N(s,v)=&\frac{\alpha_s\wedge\alpha_v}{\mid\alpha_s\wedge\alpha_v\mid}  = [-\tau_g(s)v+O(2)]T(s)\\
-&[(2{\mathcal H} (s)- {\mathcal M}(s))v+O(2)](N\wedge T)(s)
+[1+ O(2)]N(s).\endaligned $$

 \noindent Therefore it follows that  $E=<\alpha_s,\alpha_s>$, $F=<\alpha_s,\alpha_v>$, $G=<\alpha_v,
\alpha_v>$,
$e=<N,\alpha_{ss}>,$
 $\;\; f=<N,\alpha_{sv}>\;\;$ and $\; g=<N,\alpha_{vv}>$ are given by

\begin{equation}\label{eq:1f2f}
\aligned
E(s,v) &= 1-2k_g(s)v+h.o.t\\
F(s,v) &=  0+ 0.v+h.o.t\\
G(s,v) &= 1 +0.v+h.o.t\\
e(s,v)&={\mathcal M}(s) +v[\tau_g^\prime(s)-2k_g(s) {\mathcal H}(s) ]+ h.o.t\\
f(s,v) &= \tau_g(s)+ \{[2{\mathcal H}(s)-{\mathcal M}(s)]^\prime +k_g(s)\tau_g(s)\}v+ h.o.t\\
g(s,v) &=2{\mathcal H}(s)-{\mathcal M}(s) +    A(s)v+ h.o.t \endaligned\end{equation}

\begin{proposition}\label{prop:41} Let $\alpha:\mathbb V^2\to \mathbb R^3$ be an immersion of class $ C^r $, $ r\geq   6$ and $c$ be closed $ {\mathcal M}$-mean curvature line $c$  of $\alpha$,
 parametrized by arc length $s$ and  of total length $L$.
 Then the derivative of the Poincar\'e map $\pi_\alpha$ associated to $c$ is given by:

$$ln\pi_{\alpha}^\prime(0)=  \int_0^L\left[{\frac{[{\mathcal M}]_v  }{2\tau_g}} +\frac{ k_g}{\tau_g}({\mathcal H}- {\mathcal M})   \right]ds. \eqno$$

\noindent Here $\tau_g$=$ \pm \sqrt{( {\mathcal
M}-k_1)(k_2-{\mathcal M})}$.

\end{proposition}

\begin{proof} The Poincar\'e map associated to $c$ is the map $\pi_\alpha:\Sigma \to \Sigma$ defined in a transversal section to $c$ such that $\pi_{\alpha}(p)=p$ for $p\in c\cap\Sigma$ and $\pi_{\alpha}(q)$ is the first return of the  $ {\mathcal M}$-mean curvature line through $q$ to the section $\Sigma$, choosing a positive orientation for $c$. It is a local diffeomorphism and is defined, in the local chart $(s,v)$ introduced in Lemma \ref{lm:42}, by $\pi_\alpha:\{s=0\}\to \{s=L\}$, $\pi_\alpha(v_0)=v(L,v_0)$, where $v(s,v_0)$ is the solution of the Cauchy problem

$$(g-{\mathcal M})dv^2+2(f-{\mathcal M} F)dsdv+(e-{\mathcal M}E)ds^2=0, \quad v(0,v_0)=v_0.$$

Direct calculation gives
that the derivative  of the Poincar\'e map satisfies
 the following linear differential equation:

$$\frac{d}{ds}(\frac{dv}{dv_0})
 = -\frac{N_v}M(\frac{dv}{dv_0})=-\frac{ [ e-{\mathcal M}  E] _v}{  2[ f-{\mathcal M} F]  }
  (\frac{dv}{dv_0})\eqno$$

Therefore, using equation \ref{eq:1f2f} it results that
$$\frac{ [ e-{\mathcal M}  E] _v}{  2[ f-{\mathcal M}  F]  } = -\frac{\tau_g^\prime}{2\tau_g}-\frac{[{\mathcal M}  ]_v  }{2\tau_g} -\frac{k_g}{\tau_g}({\mathcal H}- {\mathcal M}). $$

Integrating  the equation above  along an arc
$[s_0,s_1]$ of  $ {\mathcal M}$-mean
curvature line,  it follows that:

 \begin{equation}\label{eq:47}
\frac{dv}{dv_0}|_{v_0=0}= \frac
{(\tau_g(s_1))^{\frac{-1}2}}{(\tau_g(s_0))^{\frac{-1}2}}exp[
\int_{s_0}^{s_1} \left[\frac{[{\mathcal M}]_v  }{2\tau_g} +\frac{ k_g}{\tau_g}({\mathcal H}- {\mathcal M})   \right]ds.\end{equation}

Applying \ref{eq:47} along the  $ {\mathcal M}$-mean curvature cycle of length $ L$, obtain
 $$\frac{dv}{dv_0}|_{v_0=0}=  exp[
  \int_0^L\left[\frac{[{\mathcal M}]_v  }{2\tau_g} +\frac{ k_g}{\tau_g}({\mathcal H}- {\mathcal M})   \right]ds. $$

>From the equation ${\mathcal K}=(eg-f^2)/(EG-F^2)$ evaluated
at $v=0$,  it follows using the expressions in   \ref {eq:1f2f}
that ${\mathcal K}={\mathcal M}[2{\mathcal H}-{\mathcal M}]-\tau_g^2.$
Developing  this equation it follows that
$\tau_g$=$ \pm \sqrt {( {\mathcal M}-k_1)(k_2-{\mathcal M})} $.

 \noindent This ends the proof. \end{proof}

 \begin{remark} \label{rem:pc} The study of the behavior of curvature lines near principal cycles was carried out in \cite{gs1}, \cite{gs2} and \cite{ag}. In this last work was established the general integral pattern for the successive derivatives of the return map.

\end{remark}

For the next theorem it is necessary to assume the additional property of  being {\it positive regular} for the function ${\mathcal M}=m(H,K)$.

\begin{definition} \label{def:bar}
A mean curvature function ${\mathcal M}=m(H,K)$, as in definition \ref{def:m}, is called {\it positive regular} if
$$ \overline{\mathcal M} ={\mathcal M}_H+2{\mathcal M}{\mathcal M}_K >0.$$
\end{definition}

\begin{proposition}\label{prop:42} Let $\alpha:\mathbb V^2\to \mathbb R^3$ be an immersion of class $ C^r $, $ r\geq   6$,  and  $c$ be a
maximal ${\mathcal M}$ -  mean curvature cycle for  $\alpha$, parametrized by arc length
and of length $L$. Consider a chart $(s,v)$ as in lemma \ref{lm:42} and consider the deformation
$$\beta_\ep(s,v)=\beta(\ep,s,v)=\alpha(s,v)+\ep [\frac{A_1(s)}6v^3]\delta (v)N(s)\eqno$$
 \noindent where $\delta =1$ in neighborhood of $v=0$, with small support and $A_1(s)=\tau_g(s)>0$.

Then $c$ is a  $ {\mathcal M}$-mean curvature cycle of $\beta_\ep$
for all $\ep $ small. Also, provided ${\mathcal M}$ is {\it
positive regular}, that is  in definition \ref{def:bar}) $
\overline{\mathcal M}>0$,  $c$ is a hyperbolic
 $ {\mathcal M}$-mean curvature cycle for $\beta_\ep$, and  $\ep\ne 0$ small.
\end{proposition}

\begin{proof}  In the chart $(s,v)$, for the immersion $\beta_\ep$, it is obtained that:

$$\aligned
E_\ep(s,v) &= 1-2k_g(s)v+h.o.t\\
F_\ep(s,v) &=  0+ 0.v+h.o.t\\
G_\ep(s,v) &= 1 +0.v+h.o.t\\
e_\ep(s,v)&={\mathcal M}(s) +v[\tau_g^\prime(s)-2k_g(s){\mathcal H}(s)\; )]+ h.o.t\\
f_\ep(s,v) &= \tau_g(s)+ [(2{\mathcal H}(s)-{\mathcal M}(s))^\prime +k_g\tau_g]v+ h.o.t\\
g_\ep(s,v) &=2{\mathcal H}(s)-{\mathcal M} (s) +v[A(s)+\epsilon A_1(s)]+ h.o.t \endaligned\eqno$$

In the expressions above $E_\ep=<\beta_s,\beta_s>$, $F_\ep=<\beta_s,\beta_v>$, $G_\ep=<\beta_v,\beta_v>$,
 \noindent $\; e_\ep=<\beta_{ss},N>$,
 $\; f_\ep=<N,\beta_{sv}>,\;$  $\; g_\ep=<N,\beta_{vv}>$, where  $N=N_\ep=\beta_s\wedge \beta_v /\mid\beta_s\wedge\beta_v\mid.$

   Let ${\mathcal M}_\ep = m({\mathcal H}_\ep,\,{\mathcal K}_\ep)$.  For all $\ep$ small it follows  that:

$$\aligned (e_\ep -{\mathcal M}_\ep  E_\ep)(s,0,\epsilon)=& 0\\
{{\mathcal K}_\ep}_v(s,0,\epsilon)=&  \ep {\mathcal M}_\ep A_1(s)+ f_1(k_g, \tau_g,  {\mathcal K} ,  {\mathcal H})(s)\\
{{\mathcal H}_\ep}_v(s,0,\epsilon)=& \frac 12  \ep  A_1(s)+ f_2(k_g, \tau_g,  {\mathcal K} ,  {\mathcal H})(s)\\
\frac{d}{d\ep }\big[{\mathcal M}_\ep]_v|_{\epsilon=0}=&\frac 12  [{\mathcal M}_H+2{\mathcal M}{\mathcal M}_K]]  A_1(s).
\endaligned
$$

Therefore  $c$ is a maximal   $ {\mathcal M}$-mean curvature cycle for all $\beta_\ep$.

Assuming that  $A_1(s)=4\tau_g(s)>0$, and also that ${\mathcal M}$ is {\it positive regular}, i.e. by definition \ref{def:bar},

$$ \overline{\mathcal M} ={\mathcal M}_H+2{\mathcal M}{\mathcal M}_K >0,$$

  \noindent it results that

$$\aligned \frac{d}{d\ep}(ln\pi^\prime(0))|_{\ep=0}=&\int_0^L
\frac{d}{d\ep}\left(\frac{({\mathcal M}_\ep)_v}{2\tau_g}+\frac{ k_g}{\tau_g}({{\mathcal H}_\ep}- {\mathcal M}_\ep)\right)ds\\
=&\int_0^L \tau_g{\overline{\mathcal M}}ds > 0.  \endaligned $$
\end{proof}

As a synthesis of propositions \ref{prop:41} and \ref{prop:42},  the following theorem is obtained.

\begin{theorem}\label{th:41} Let $\mathcal M$ be a positive regular mean curvature function.
An immersion $\alpha\in {\mathcal I}^{r,s}({\mathbb V}^2)$, $ r\geq   6$, is  $C^6-$local  $ {\mathcal M}$-mean curvature structurally stable at a  $ {\mathcal M}$-mean curvature cycle $c$ if only if,

$$  \int_0^L\left[\frac{[{\mathcal M}]_v  }{2\tau_g} +\frac{ k_g}{\tau_g}({\mathcal H}- {\mathcal M})   \right]ds \neq 0. $$
\end{theorem}

\begin{proof} Using propositions  \ref{prop:41} and \ref{prop:42}, the  local topological character of the foliation can be changed by small perturbation of the immersion, when the cycle is not hyperbolic.
\end{proof}

\section{ $ {\mathcal M}$-Mean Curvature Lines near the Parabolic Curve}\label{sec:5}

In this section will be studied the behavior of the $ {\mathcal M}$-mean curvature lines near the parabolic points of an immersion, assuming that the quadratic differential equation \ref{eq:harm} is univalued there. This is done under  two regularity conditions imposed  in definitions \ref{def:mb1} and  \ref{def:mb2}.  The  motivation comes from the previous study  of the classical harmonic and geometric mean curvature functions; see \cite {h, m}.

\begin{definition}\label{def:mb1}
A mean curvature function $\mathcal M=m(H,K)$ is called {\it 1-regular}   if either
\begin{itemize}
\item[a)] $m(H,0)=0$ and $({\partial m}/{\partial K})(H,0)> 1/(2H)>0$,
or
\item[b)] $m(H,0)=2H$ and $({\partial m}/{\partial K})(H,0)<- 1/(2H)<0$.
\end{itemize}
\end{definition}

 \begin{remark}\label{rem:3}
For mean curvature functions $m$, with $m(H,0)=0$, it always holds that  $({\partial m}/{\partial K})(H,0)\geq 1/(2H)>0.$
\noindent The $1-$ regular condition states that the inequality is strict.
In fact,
$$\frac{\partial m}{\partial K}(H,0)=\lim_{K\to 0}\frac{m(H,K)-M(H,0)}{K}\geq \lim_{K\to 0} \frac{H-\sqrt{H^2-K}}{K}=\frac{1}{2H}. $$

Analogously for the   case where $m(H,0)=2H$.
\end{remark}

\begin{definition}\label{def:mb2}
 A mean curvature function $\mathcal M=m(H,K)$ is called {\it 1/2-regular}  if either
\begin{itemize}
\item[a)]  $m(H,0)=0$, $m(H,K)=\overline{m}(H,\sqrt{K})$ for some analytic function $\overline{m}(H,S)$ which furthermore satisfies
$({\partial\overline{m}}/{\partial S})(H,0)>0$,  or

\item[b)]
$m(H,0)=2H$, $m(H,K)=\overline{m}(H,\sqrt{K})$ for some analytic
function $\overline{m}(H,S)$ which furthermore satisfies
$({\partial \overline{ m}}/{\partial S})(H,0)<0$.

\end{itemize}
\end{definition}

The natural examples for cases a) in the definitions above  are the Harmonic( $m=K/H$)  and Geometric ($m=\sqrt{K}$) mean curvatures.
For cases b), take $m= 2H-K/H$ and $m= 2H-\sqrt{K}$.

\subsection{ $ {\mathcal M}$-mean curvature lines near a parabolic line: the 1-regular,  case a)}

 Let  $0$ be a parabolic  point of a $ C^r , \; r\geq   6$, immersion  $\alpha$ parametrized in a Monge chart $(x,y)$ by $\alpha(x,y)=(x,y,z(x,y))$, where

\begin{equation}\label{eq:ez}
\aligned z(x,y) =&\frac k2 y^2 + \frac a6 x^3+\frac b2 xy^2+\frac d2 x^2y+\frac c6 y^3 + \frac A{24}x^4+\frac B{6} x^3y\\
+&\frac{C}4 x^2y^2+\frac D6 xy^3+\frac E{24}y^4+\sum_{i+j=5} r_{ij} \frac{x^iy^j }{i! j!}+O(6)\endaligned
\end{equation}

The coefficients of the first and second fundamental forms are given by:

\begin{equation}\label{eq:1f2fp}
\aligned E(x,y)=& 1+O(4)\\
F(x,y)=&  \frac 12 ak x^2y+ kd xy^2+\frac kb y^3+O(4)\\
G(x,y)=&1+k^2 y^2 +kd x^2y +2kb xy^2+kc y^3+0(4)\\
e(x,y)=& ax+dy+\frac 12 A x^2+Bxy+\frac 12 Cy^2 +\frac 16 r_{50} x^3
 +\frac 12 r_{41}x^2y\\
+& \frac 12 (r_{32}-ak^2)xy^2
+\frac 16 (r_{23}-3k^2d)y^3+O(4)\\
f(x,y)=& dx+ by+\frac 12 Bx^2+Cxy+\frac 12 D y^2 +\frac 16 r_{41} x^3
 +  \frac 12 r_{32}x^2y\\
+& \frac 12 (r_{23}- k^2d)xy^2
+\frac 16 (r_{14}-3k^2b)y^3+O(4) \\
g(x,y)=& k+bx+cy+\frac 12 Cx^2+Dxy+\frac 12{( E-k^3)} y^2 +\frac 16 r_{32} x^3\\
 +& \frac 12 (r_{23}-k^2d)x^2y+ \frac 12 (r_{14}-3b k^2)xy^2\\
+&\frac 16 (r_{05}-6k^2c)y^3+O(4)\endaligned
\end{equation}

The Gaussian   and the Mean curvatures are given by

\begin{equation}\label{eq:gc}
 \aligned
{\mathcal K}(x,y)=& k(ax+dy)+ \frac 12(Ak+2ab-2d^2)x^2+(Bk+ac-bd)xy\\
+&\frac 12(Ck+2cd-2b^2)y^2+ \frac 16( kr_{50} +3Ab+3aC-6Bd)x^3\\
+& \frac 12 ( kr_{32}-4ak^3+2cB+aE-3bC)xy^2
\\
+&\frac 12(kr_{41}+cA+2aD-3Cd)x^2y\\
+& \frac 16(kr_{23}+3cC+3Ed-6bD-12k^3d)y^3+O(4),\\
{\mathcal H}(x,y)=& \frac 12k+\frac 12(a+b)x+\frac 12 (c+d)y+\frac 14(A+C ) {x^2} \\
+& \frac 12 (B+D)xy +\frac 14( E +C-3 k^3)  {y^2}
+\frac 1{12}(r_{50}+r_{32})x^3\\
 +&\frac 14 ( r_{32}+r_{14}+k^2(2d-9b-a))xy^2\\
+&\frac 14 (r_{41}+r_{23}+k^2(a-3d))x^2y\\
+&\frac 1{12} ( r_{23}+r_{05}+k^2(3b-9c-3d))y^3+O(4)
\endaligned \end{equation}

Let $\mathcal M$ be a 1- regular function. Write
$${\mathcal M}(x,y)= m({\mathcal H},{\mathcal K})(x,y)= (m_0+m_1x+m_2y+O(2)){\mathcal K}(x,y).$$

The coefficients of the quadratic differential equation (\ref{eq:harm}) are given by $L,\, M,\, N$ as follows:

$$L=  g-{\mathcal M}G,  \;\; M= 2(f-{\mathcal M}F), \;\; N=  e-{\mathcal M}E.$$

\begin{equation}\label{eq:lmn}
\aligned
L=& k +(b-akm_0)x+(c-km_0d)y \\
+&\frac 12 [ C+(2d^2-2ab- Ak)m_0-2akm_1]{x^2} \\
+&[ D+(bd-ac-kB)m_0-km_1d-akm_2]xy \\
+&\frac 12[(2b^2- kC-2cd)m_0- k^3+ E-2km_2d] {y^2}  +O(3)\\
M =& 2d x+2b y +B {x^2} +2Cxy +  D  y^2 +O(3)\\
 N=&a (1-km_0)x+ d(1-km_0)y \\
+&\frac 12[(2d^2 -2ab- Ak)m_0+ A -2akm_1] {x^2}  \\
+& [B +(bd-ac-kB)m_0-km_1d-akm_2]xy \\
+&\frac 12[(2b^2-2cd-
kC)m_0-km_2d+ C] {y^2} +O(3)\\
\endaligned
\end{equation}

\begin{lemma}\label{lm:p1} Let $0$ be a parabolic point and consider the parametrization $(x,y,z(x,y))$ as  above. If $k>0$ and $a^2 + d^2 \ne 0$ then the set of parabolic points is locally a regular curve normal to the vector $(a,d)$ at $0$.
\begin{itemize}
\item[i)]
 If $a\ne 0$ the parabolic curve is transversal to the minimal principal direction $(1,0)$.
\item[ii)]
If $a=0$ then the parabolic curve is tangent to the principal direction given by $(1,0)$ and has quadratic contact with the corresponding minimal principal curvature line if $dk(Ak-3d^2)\ne 0$.
\end{itemize}
\end{lemma}

\begin{proof} If $a \ne 0$, from the expression of $\mathcal K$ given by equation \ref{eq:gc} it follows that the parabolic line is given by $x=-\frac{d}{a}y+O_1(2)$ and so is transversal to the principal direction $(1,0)$ at $(0,0)$.

 If $a=0$, from the expression of $\mathcal K$ given by equation \ref{eq:gc} it follows that the parabolic line is given by $y=\frac{ 2d^2-Ak}{2dk}x^2+O_2(3)$  and that $y=-\frac{d}{2k}x^2+O_3(3)$ is the principal line tangent to the principal direction $(1,0)$.
Now the condition of quadratic contact $\frac{ 2d^2-Ak}{2dk}\ne -\frac{d}{2k}$ is equivalent to $dk(Ak-3d^2)\ne 0$.\end{proof}

\begin{figure}[htbp]
 \begin{center}
 \includegraphics[angle=0, width=11cm]{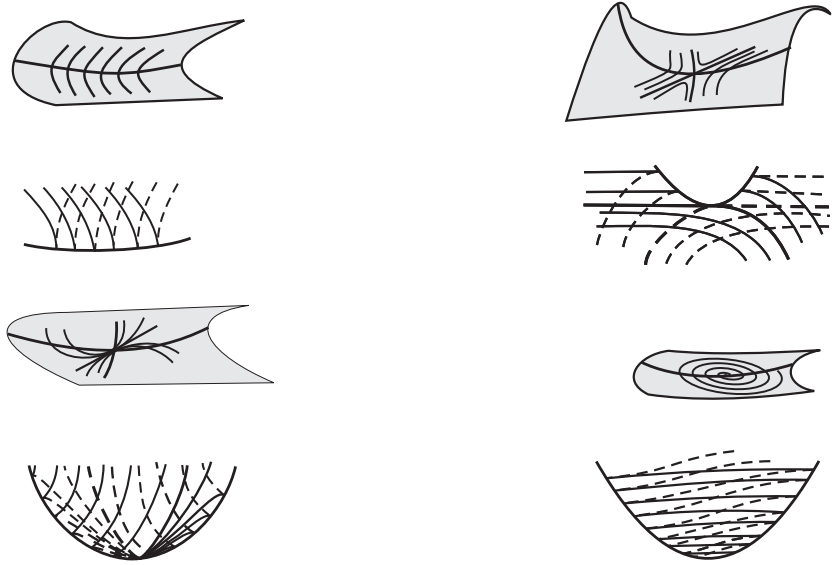}
 \caption{$ {\mathcal M}$-mean curvature lines near a parabolic point (cuspidal, folded saddle,    folded node and folded focus) and
  their separatrices \label{fig:mparab} }
    \end{center}
  \end{figure}

\begin{proposition}\label{prop:p1} Let $0$ be a parabolic point and the Monge chart $(x,y)$ as above.
Suppose $\mathcal M$ is $1-$regular at $(k,0)$
with ${\partial m}/{\partial K}(k,0)=m_0>1/k$.
\begin{itemize}
\item[i)]
If $a\ne 0$ then the mean   $ {\mathcal M}$-curvature lines are transversal to the parabolic curve and the mean curvatures lines are shown in the Figure \ref{fig:mparab}, the cuspidal case.

\item[ii)]
If $a=0$ and $\sigma=  (A k-3 d ^2)\ne 0$ then the  mean   $\mathcal M$-
curvature lines are shown in the Figure \ref{fig:mparab}. In fact, if
$\sigma >0$ then the $ {\mathcal M}$-mean    curvature lines are folded
saddles. Otherwise,  if $\sigma <0$ then the $ {\mathcal M}$-mean   curvature lines are folded nodes or folded focus according to
$\delta=[d^2(km_0-25)+8Ak]$
 be positive or negative. The two separatrices of these
tangential singularities, folded saddle and folded node, as
illustrated in the Figure \ref{fig:mparab}
 below, are called parabolic
separatrices.
\end{itemize}
\end{proposition}

\begin{proof}
Consider the quadratic differential equation
$$H(x,y,[dx:dy])=Ldy^2+Mdxdy +Ndx^2 =0$$

\noindent and the Lie-Cartan line field $X$ of class $C ^{r-3}$ defined by

$$\aligned x^\prime =&H_p\\
y^\prime =& pH_p\\
p^\prime =&-(H_x+pH_y), \hskip 1cm p=\frac{dy}{dx}\endaligned$$
\noindent where $L$, $M$ and $N$ are given by equation \ref{eq:lmn}.

If $a\ne 0$ the vector $Y$ is regular and therefore  the $ {\mathcal M}$-mean    curvature lines are transversal to the parabolic line and at parabolic points these lines are  tangent to the  principal direction $(1,0)$.

If $a=0$, direct calculation  gives $H(0)=0,\;\; H_x(0)=0, \;\; H_y(0)=-kd, \;\; H_p(0)=0.$

 \begin{equation}\aligned
DX(0)=\left(\begin{matrix} 2d & 2b & 2k\\  0 & 0 & 0\\a_{31}&a_{32}&a_{33}\end{matrix}\right)
\endaligned
\end{equation}
\noindent where,

 $$a_{31}=(Ak-2d^2)m_0-A,\; a_{32}=(kB-bd)m_0-B+kdm_1+,\; a_{33}=(km_0-3)d.$$

The non vanishing eigenvalues of $DX(0)$ are

$$\lambda_1=\frac{1}2\big[\; d(km_0-1)-\sqrt{(-1+km_0)(d^2(km_0-25)+8Ak)}\;\big],$$

$$\lambda_2= \frac{1}2\big[\; d(km_0-1)+\sqrt{(-1+km_0)(d^2(km_0-25)+8Ak)}\;\big]$$

Therefore, $\lambda_1\lambda_2=2(1-km_0)(Ak-3d^2)$.

It follows that $0$ is a hyperbolic  singularity provided $\sigma=
(Ak-3d^2)\ne 0$. If $\sigma >0$ then the $ {\mathcal M}$-mean
curvature lines are folded saddles and   if $\sigma <0$ then the
$ {\mathcal M}$-mean curvature lines are folded nodes
$(\;[d^2(km_0-25)+8Ak]>0 )$ or folded focus
$(\;[d^2(km_0-25)+8Ak]>0 )$. See Figure  \ref{fig:mparab}.
 \end{proof}

\begin{theorem}\label{th:51a}
Assume that the $ {\mathcal M}$-mean curvature function $\mathcal
M$ is $1$-regular as in definition \ref{def:mb1}, case a). An immersion
$\alpha\in {\mathcal I}^{r,s}({\mathbb V}^2)$, $ r\geq   6$, is
$C^6-$local  $ {\mathcal M}$-mean curvature structurally stable at
a tangential parabolic point $p$  if only if, the condition
$\sigma\delta \neq 0$ in  proposition \ref{prop:p1}  holds.

\end{theorem}

\begin{proof}
Direct from Lemma   \ref{lm:p1} and proposition \ref{prop:p1}, the  local topological character of the foliation can be changed by small perturbation of the immersion when $\delta \sigma = 0$.
\end{proof}

\subsection{ $ {\mathcal M}$-mean curvature lines near a parabolic line: the  1/2-regular, case a)}\label{sb:52}\hskip 8cm
\newline

 Let  $0$ be a parabolic  point of a $ C^r , \; r\geq   6$, immersion  $\alpha$ parametrized in a Monge chart $(x,y)$ by $\alpha(x,y)=(x,y,z(x,y))$, where $z $ is as in equation \ref{eq:ez}.

The coefficients of the first and second fundamental forms are given by expressions \ref{eq:1f2fp}.

The Gaussian and Arithmetic Mean curvatures are given by equations
\ref{eq:gc}.

Below is established the typical behavior of $\mathcal M$-
mean curvature lines for a function ${\mathcal M}={\mathcal M_1}\sqrt{\mathcal K} $,  as in definition \ref{def:mb2}.

Squaring both members of  the differential equation $k_n(x,y,[dx:dy]) ={\mathcal M_1}\sqrt{\mathcal K} $ to remove the square root singularity,  gives the following
 quartic differential equation:

\begin{equation}\label{eq:mgc4} A_{40} dx^4+ A_{31}dx^3dy + A_{22}dx^2dy^2+A_{13} dxdy^3+A_{04} dy^4 =0
\end{equation}

\noindent where,

$$\aligned A_{40}&= e^2 (E G-F^2)-E^2 (eg-f^2){\mathcal M}_1^2\\
A_{31}&=4ef(EG-F^2)-4EF(eg-f^2){\mathcal M}_1^2\\
A_{22}&=( 4f^2+2eg)(EG-F^2)-(2EG+4F^2)(eg-f^2){\mathcal M}_1^2\\
A_{13}&=4fg(EG- F^2)-4FG(eg-f^2){\mathcal M}_1^2 \\
A_{ 04}&= g^2 (E G-F^2)-G^2 (eg-f^2){\mathcal M}_1^2 \endaligned$$

Writing
$${\mathcal M}_1(x,y)=m({\mathcal H}(x,y), {\mathcal K}(x,y))=m_0+m_1x+m_2y+O(2),$$

\noindent the coefficients of the quartic differential equation (\ref{eq:mgc4}) are given by

\begin{equation}\label{eq:ai}
\aligned A_{40}=&-km_0^2(ax+dy)+\frac 12[2a^2+m_0^2(2d^2-2ab-Ak)-4akm_0m_1]x^2\\
+& [2ad-(Bk+bd-ac)m_0^2-2km_0(dm_1+am_2)]xy\\
+& \frac 12[2d^2+m_0^2(2b^2  -2cd-Ck)-4kdm_0m_2]y^2+0(3)\\
%
A_{31} =& 4 da x^2+4( d^2+ab)xy+4 d by^2+ 0(3)\\
%
A_{22}=& 2 k(1-m_0^2)(ax+dy)\\
+&[ m_0^2(2d^2-2ab-kA)-4akm_0m_1+Ak+2ab+4d^2]x^2\\
+&[10bd+2ac+2kB-4km_0(m_1d+am_2)+2m_0^2(bd-kB-ac)]xy\\
+&[(2b^2-kC-2cd)m_0^2 -4m_0m_2kd+kC+2cd+4b^2]y^2+0(3)\\
%
A_{13}=& 4k(dx+by)+ (2Bk+4bd)x^2+4(Ck+cd+b^2)xy\\
+&(2kD+4bc)y^2+0(3)\\
%
A_{04}=& k^2+k(2b-am_0^2)x+k(2c-dm_0^2)y\\
+&\frac 12[2b^2+2kC-4akm_0m_1+ m_0^2(2d^2-2ab -Ak)]x^2\\
+&[2kD+2bc+m_0^2(bd-kB-ac)-2km_0(m_1d+m_2a)]xy\\
+&\frac 12[2c^2-2k^4+  2kE +m_0^2(2b^2-kC-2cd)-4m_0m_2 d ]y^2+0(3)\\
\endaligned
\end{equation}

\begin{lemma}\label{lem:p5} Let $0$ be a parabolic point and consider the parametrization $(x,y,z(x,y))$ as  above. If $k>0$ and $a^2 + d^2 \ne 0$ then the set of parabolic points is locally a regular curve normal to the vector $(a,d)$ at $0$.
\begin{itemize}
\item[i)]
 If $a\ne 0$ the parabolic curve is transversal to the minimal principal direction $(1,0)$.

\item[ii)] If $a=0$ then the parabolic curve is tangent to the principal direction given by $(1,0)$ and has quadratic contact with the corresponding minimal principal curvature line if $dk(Ak-3d^2)m_0^2\ne 0$.

\end{itemize}
\end{lemma}

\begin{proof} If $a \ne 0$, from the expression of $\mathcal K$ given by equation \ref{eq:gc} it follows that the parabolic line is given by $x=-\frac{d}{a}y+O_1(2)$ and so is transversal to the principal direction $(1,0)$ at $(0,0)$.

 If $a=0$, from the expression of $\mathcal K$ given by equation
\ref{eq:gc} it follows that the parabolic line is given by
 $y=\frac{ 2d^2-Ak}{2dk}x^2+O_2(3)$  and that $y=-\frac{d}{2k}x^2+O_3(3)$ is the principal line tangent to the principal direction $(1,0)$.
Now the condition of quadratic contact $\frac{ 2d^2-Ak}{2dk}\ne -\frac{d}{2k}$ is equivalent to $dk(Ak-3d^2)\ne 0$.\end{proof}

\begin{figure}[htbp]
\begin{center}
\includegraphics[angle=0, width=11cm]{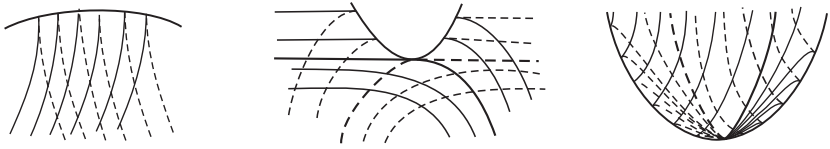}
       \caption {$ {\mathcal M}$-mean curvature lines near a parabolic point
       (cuspidal, folded saddle and folded node) and their separatrices \label{fig:mpr}}
   \end{center}
 \end{figure}

\begin{proposition}\label{prop:p2} Let $0$ be a parabolic point and the Monge chart $(x,y)$ as
above and   $\mathcal M $  be a mean curvature function {\it
1/2-regular}, case a).
\begin{itemize}
\item[i)]
If $a\ne 0$ then the $ {\mathcal M}$-mean  curvature lines are
transversal to the parabolic curve, as  shown in  Figure
\ref{fig:mpr}, the cuspidal case.

\item[ii)]
If $a=0$ and $\sigma=dk(Ak-3d^2)m_0^2\ne 0$ then the  $\mathcal M$- mean
curvature lines are shown in the Figure \ref{fig:mpr}. In fact, if
$\sigma >0$ then the $\mathcal M$- mean   curvature lines are folded
saddles. Otherwise,  if $\sigma <0$ then the $ {\mathcal M}$-mean
curvature lines are folded nodes. The two separatrices of these
tangential singularities, folded saddle and folded node, as
illustrated in Figure \ref{fig:mpr}
are called parabolic
separatrices.
\end{itemize}
\end{proposition}

\begin{proof}
Consider the quartic differential equation,

$$H(x,y,\frac{dy}{dx})=A_{04}(\frac{dy}{dx})^4+A_{13}(\frac{dy}{dx})^3+A_{22}(\frac{dy}{dx})^2+A_{31}\frac{dy}{dx}+A_{40}=0$$

\noindent and the Lie-Cartan line field of class $C^{r-3}$ defined by

$$\aligned x^\prime =&H_p\\
y^\prime =& pH_p\\
p^\prime =&-(H_x+pH_y), \hskip 1cm p=\frac{dy}{dx}\endaligned$$
\noindent where $A_{ij}$ are given by equation \ref{eq:ai}.

If $a\ne 0$ the vector $Y$ is regular and therefore  the  $ {\mathcal M}$-mean curvature lines are transversal to the parabolic curve.
 If $a = 0$, the parabolic curve is tangent to the  principal direction $(1,0)$.

For $a=0$, direct calculation  gives $H(0)=0,\;\; H_x(0)=0, \;\; H_y(0)=-kd, \;\; H_p(0)=0.$

Therefore, solving the equation $H(x,y(x,p),p)=0$ near $0$  it follows,  by the Implicit Function Theorem, that:
$$y=y(x,p)= \frac{2d^2-Ak}{2kd} x^2-
\frac{  r_{50}k^2d+6Akbd-3BAk^2-6d^3b }{6k^2d^2}x^3+O(4).$$

Therefore the vector field
$Y$ given by  the differential equation below
$$\aligned x^\prime =&H_p(x,y(x,p),p)\\
p^\prime =&-(H_x+pH_y)(x,y(x,p),p)\endaligned$$
\noindent is given by
$$\aligned x^\prime =&  \frac{4d^3}k  x^3+12d^2 x^2p+12kd xp^2+4k^2p^3+O(4)\\
p^\prime =& (Ak-2d^2)m_0^2x+kdm_0^2p+O(2).\endaligned
$$

The singular point $0$ is  isolated and the eigenvalues of the
linear part of $Y$ are given by $\lambda_1=0$ and $\lambda_2=m_0^2kd$.
The correspondent eigenvectors are given by $f_1=(1,(2d^2-Ak)/dk)$
and $f_2=(0,1)$.

Performing the calculations, restricting $Y$ to the center manifold $W^c$ of class $C ^{r-3}$, $T_0W^c=f_1$,  it follows that
$$Y_{c}= -\frac 23 \frac{(Ak-3d^2)^3}{kd^3}x^3 +0(4)$$

It follows that $0$ is a topological saddle or node  of cubic type
provided $\sigma (Ak-3d^2)km_0^2d\ne 0$. If $\sigma >0$ then the
$\mathcal M$- mean   curvature lines are folded saddles and   if
$\sigma <0$ then the   $\mathcal M$- mean   curvature lines are
folded nodes. In the case $\sigma>0$, the center manifold $W^c$ is
unique, \cite{sac}, cap. $V$, page 319, and so the saddle
separatrices      are well defined. See Figure \ref{fig:selano}
 below.

\begin{figure}[htbp]
\begin{center}
\includegraphics[angle=0, width=10cm]{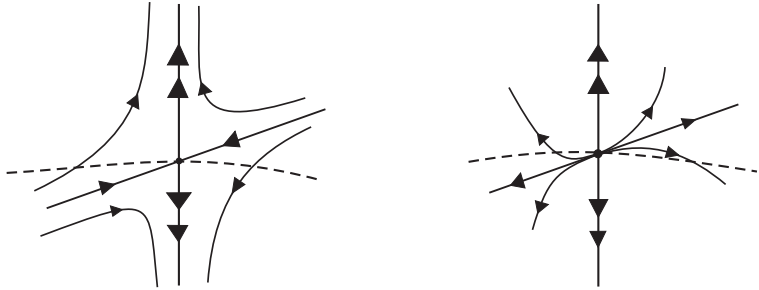}
  \caption {Phase portrait of the vector field $Y$ near singularities \label{fig:selano} }
   \end{center}
 \end{figure}

Notice that due to the constrains of the problem treated here, the non
  hyperbolic saddles and nodes, which in the standard theory would
  bifurcate into three singularities, are actually structurally stable (do
  not bifurcate).

 \end{proof}

\begin{remark}

The reader may find a more complete study of the partial hyperbolicity  structure in the theorem above, which can be expressed in the context of normal hyperbolicity,  in the paper of Palis and Takens \cite {pt}.

 For a deeper analysis of the lost of the hyperbolicity condition and the
  consequent bifurcations,  the reader is addressed to the book of Roussarie
 \cite{rr}.

\end{remark}

\begin{remark} As in the analysis of geometric mean curvature lines near a parabolic point, see \cite{g}, the singularity of the Lie-Cartan vector field above is semi-hyperbolic of order $3$. However,  the terms of third order of the functions $A_{ij}$ of equation (\ref{eq:ai}) have no contribution to  the orbit structure  around  the singularity. This was  confirmed by computer algebraic  calculations.
\end{remark}

\begin{theorem}\label{th:51}
Assume that the $ {\mathcal M}$-mean curvature function $\mathcal M$ is $1/2$-regular as in definition \ref{def:mb2}, case a).

An immersion $\alpha\in {\mathcal I}^{r,s}({\mathbb V^2)}$, $ r\geq   6$, is  $C^6-$local $ {\mathcal M}$- mean curvature structurally stable at a tangential parabolic point $p$
if only if, the condition $\sigma \neq 0$ in  proposition \ref{prop:p2}  holds.

\end{theorem}

\begin{proof}
Direct from Lemma   \ref{lem:p5} and proposition \ref{prop:p2}, the  local topological character of the foliation can be changed by small perturbation of the immersion when $\sigma = 0$.
\end{proof}

\subsection{ $ {\mathcal M}$- mean curvature lines near a parabolic line: the 1-regular, case b)}\hskip 8cm
\newline

 Let  $0$ be a parabolic  point of a $ C^r , \; r\geq   6$, immersion  $\alpha$ parametrized in a Monge chart $(x,y)$ by $\alpha(x,y)=(x,y,z(x,y))$, where

\begin{equation}\label{eq:ez2}\aligned z(x,y) =&\frac k2 x^2 + \frac a6 x^3+\frac b2 xy^2+\frac d2 x^2y+\frac c6 y^3 +\frac A{24}x^4+\frac B{6} x^3y\\
+&\frac{C}4 x^2y^2+\frac D6 xy^3+\frac E{24}y^4+\sum_{i+j=5} r_{ij} \frac{x^iy^j }{i! j!}+O(6)\endaligned
\end{equation}

The coefficients of the first and second fundamental forms are given by:

\begin{equation}\label{eq:1f2fp2}\aligned E(x,y)=&
1+k^2x^2+kax^3+2kd x^2y+kbxy^2+O(4)\\
F(x,y)=& \frac 12 kdx^3+kbx^2y+\frac 12 kcxy^2 +
O(3)\\
G(x,y)=&1+ O(4)\\
e(x,y)=& k+ax+dy+\frac 12 (A-k^3)x^2+Bxy+\frac 12 Cy^2\\
+&\frac 16(r_{50}-6ak^2){x^3}+ (r_{41}-3k^2d)x^2y\\
+&\frac 12(r_{32}-k^2b){xy^2}+\frac 16 r_{23}y^3+O(4)\\
f(x,y)=& dx+by+\frac 12 Bx^2+Cxy+\frac 12 Dy^2\\
+& \frac 16(r_{41}-3dk^2)x^3+(r_{32}-k^2b)x^2y\\
+&\frac 12 r_{23}xy^2+\frac 16 r_{14}y^3+O(4)\\
g(x,y)=& bx+cy+\frac 12 Cx^2+Dxy+\frac 12 Ey^2\\
+&\frac 16(r_{32}-3k^2b)x^3+\frac 12(r_{23}-k^2c)x^2y\\
+&\frac 12 r_{14}xy^2+\frac 16 r_{05}y^3+O(4)\endaligned
\end{equation}

The Gaussian   and the Mean curvatures are given by

\begin{equation}\label{eq:gc2}
 \aligned
{\mathcal K}(x,y)=& k(bx+cy)+\frac 12(2ab+kC-2d^2)x^2+(ac-bd+kD)xy\\
+&\frac 12(2cd+kE-2b^2)y^2 +\frac 16(3Ab-12k^3b+3aC+kr_{32}-6dB)x^3\\
+&\frac 12(2aD-4k^3c-3dC+Ac+r_{23}k)x^2y\\
+&\frac 12(2Bc+aE-3bC +r_{14}k)xy^2\\
+&\frac 16( 3cC-6bD+3dE+kr_{05})y^3+O(4),\\
{\mathcal H}(x,y)=& \frac 12 k+\frac 12(b+a)x+\frac 12(c+d)y+\frac 14(C+A-3k^3)x^2\\
+&\frac 12(D+B)xy
+\frac 14(E+C)y^2
+\frac{ 1}{12}(r_{50}+r_{32}-3k^2(b+6a)){x^3}\\
+ &\frac 14(r_{41}+r_{23}-k^2(c+9d))x^2y
+\frac 14(r_{32}+r_{14}-3k^2b)xy^2\\
+& \frac 1{12}( r_{23}+r_{05})y^3+O(4)
\endaligned \end{equation}

Let $\mathcal M$ be a 1- regular function. Write
$${\mathcal M}(x,y)= m({\mathcal H},{\mathcal K})(x,y)= 2{\mathcal H}-(m_0+m_1x+m_2y+O(2)){\mathcal K}(x,y).$$

The coefficients of the quadratic differential equation (\ref{eq:harm}) are given by $L,\, M,\, N$ as follows:

$$L=  g-{\mathcal M}G,  \;\; M= 2(f-{\mathcal M}F), \;\; N=  e-{\mathcal M}E.$$

\begin{equation}\label{eq:lmn2}
\aligned
 L=&-k+ (kbm_0-a)x+(kcm_0-d)y\\
+& \frac 12 [m_0({kC}+2ab-2d^2)+2kbm_1+   {3k^3}-A]{x^2}\\
+&[m_0(ac+kD-bd) +k(cm_1+bm_2)-B]xy\\
+&\frac 12[m_0({kE}+2cd-2b^2) -C+2kcm_2]{y^2}+O(3)\\
M=& 2dx+2by+B  {x^2} +2Cxy+D  {y^2}+O(3)\\
N=& (-1+km_0)(bx+cy)+ \frac 12 [m_0(2ab+{ kC}-2d^2) +2kbm_1-C]{x^2}\\
+& [m_0(ac+kD-bd)  +k(bm_2+cm_1)-D]xy\\
+& \frac 12[m_0(2cd-2b^2+{ k  E}) +2k c m_2-E]{y^2}+O(3)
\endaligned
\end{equation}

\begin{lemma}\label{lem:p1} Let $0$ be a parabolic point and consider the parametrization $(x,y,z(x,y))$ as  above. If $k>0$ and $b^2 + c^2 \ne 0$ then the set of parabolic points is locally a regular curve normal to the vector $(b,c)$ at $0$.
\begin{itemize}
\item[i)]
 If $b\ne 0$ the parabolic curve is transversal to the maximal principal direction $(1,0)$.
\item[ii)]
If $b=0$ then the parabolic curve is tangent to the maximal principal direction given by $(1,0)$ and has quadratic contact with the corresponding maximal principal curvature line if $  d(2d-c)-kC \ne 0$.
\end{itemize}
\end{lemma}

\begin{proof} If $b\ne 0$, from the expression of $\mathcal K$ given by equation \ref{eq:gc} it follows that the parabolic line is given by $y=-\frac{b}{c}x+O_1(2)$ and so is transversal to the maximal principal direction $(1,0)$ at $(0,0)$.

 If $b=0$, from the expression of $\mathcal K$ given by equation \ref{eq:gc} it follows that the parabolic line is given by $y=\frac{ 2d^2-Ck}{2ck}x^2+O_2(3)$  and that $y=\frac{d}{2k}x^2+O_3(3)$ is the principal line tangent to the principal direction $(1,0)$.
Now the condition of quadratic contact $\frac{ 2d^2-kC}{2ck}\ne \frac{d}{2k}$ is equivalent to $ d(2d-c)-kC \ne 0$.\end{proof}

\begin{proposition}\label{prop:p1b} Let $0$ be a parabolic point and the Monge chart $(x,y)$ as above.
Suppose  $\mathcal M$ is $1-$regular, case b), at $(k,0)$
with ${\partial m}/{\partial K}(k,0)=-m_0< -1/k<0$.
\begin{itemize}

\item[i)]
If $b\ne 0$ then the $\mathcal M$- mean    curvature lines are transversal to the parabolic curve and the mean curvatures lines are shown in the Figure \ref{fig:mparab}, the cuspidal case.

\item[ii)]
If $b=0$ and  $\sigma= (cd-2d^2+kC)kc m_0^2\ne 0$ then the $\mathcal M$- mean
curvature lines are shown in the Figure \ref{fig:mparab}. In fact, if
$\sigma >0$ then the $\mathcal M$-mean    curvature lines are folded
saddles. Otherwise,  if $\sigma <0$ then the mean   $\mathcal M$-
curvature lines are folded nodes or folded focus according to
$\delta= (k m_0-1)(c^2+8kC +8d(c-2d))$
 be positive or negative. The two separatrices of these
tangential singularities, folded saddle and folded node, as
illustrated in the Figure \ref{fig:mparab}
 of proposition \ref{prop:p1}, are called parabolic
separatrices.
\end{itemize}
\end{proposition}

\begin{proof}
Consider the quadratic differential equation
$$H(x,y,[dx:dy])=Ldy^2+Mdxdy +Ndx^2 =0$$

\noindent and the Lie-Cartan line field $X$ of class $C ^{r-3}$ defined by

$$\aligned x^\prime =&H_p\\
y^\prime =& pH_p\\
p^\prime =&-(H_x+pH_y), \hskip 1cm p=\frac{dy}{dx}\endaligned$$

\noindent where $L$, $M$ and $N$ are given by equation \ref{eq:lmn2}.

If $b\ne 0$ the vector $Y$ is regular and therefore  the $\mathcal M$- mean
 curvature lines are transversal to the parabolic line and at parabolic points
 these lines are  tangent to the  principal direction $(1,0)$.

If $b=0$, direct calculation  gives $H(0)=0,\;\; H_x(0)=0, \;\; H_y(0)=c(km_0-1), \;\; H_p(0)=0.$

 \begin{equation}\aligned
DX(0)=\left(\begin{matrix} 2d & 0 & -2k\\  0 & 0 & 0\\a_{31}&a_{32}&a_{33}\end{matrix}\right)
\endaligned
\end{equation}
\noindent where,

 $$\aligned a_{31} =& C+m_0(2d^2-kC), \;\;\;\; a_{32}=D-(kD+ac)m_0-kcm_1,\\
a_{33}=&c-2d-kcm_0. \endaligned $$

The non vanishing eigenvalues of $DX(0)$ are

$$\lambda_1=\frac{1}2\big[\;  c(1-km_0)-\sqrt{(-1+km_0)(km_0c^2-c^2+8kC +8d(c-2d)) }\;\big],$$

$$\lambda_2= \frac{1}2\big[\; c(1-km_0)+\sqrt{(-1+km_0)(km_0c^2-c^2+8kC +8d(c-2d))}\;\big]$$

Therefore, $\lambda_1\lambda_2= 2 ( 1-km_0)(kC+cd-2d^2)$.

It follows that $0$ is a hyperbolic  singularity provided $\sigma= (cd-2d^2+kC)kc m_0^2\ne 0$. If $\sigma >0$ then the $\mathcal M$-mean  curvature lines are folded saddles and   if $\sigma <0$ then the      $\mathcal M$-mean curvature lines are folded nodes $(\;(c^2+8kC +8d(c-2d))>0 )$ or folded focus $(\;(c^2+8kC +8d(c-2d)) <0 )$.
 \end{proof}

\begin{theorem}\label{th:51ab}
Assume that the  $\mathcal M$-mean curvature function $\mathcal M$ is $1$-regular as in definition \ref{def:mb1}, case b).
An immersion $\alpha\in {\mathcal I}^{r,s}({\mathbb V}^2)$, $ r\geq   6$, is  $C^6-$local   $\mathcal M$-mean curvature structurally stable at a tangential parabolic point $p$
if only if, the condition $\sigma\delta \neq 0$ in  proposition \ref{prop:p1}  holds.

\end{theorem}

\begin{proof}
Direct from Lemma   \ref{lem:p5} and proposition \ref{prop:p1b}, the  local topological character of the foliation can be changed by small perturbations of the immersion when $\delta \sigma = 0$.
\end{proof}

\subsection{$ {\mathcal M}$- mean curvature lines near a parabolic line: the 1/2-regular,  case b)}\hskip 8cm
\newline

Below are formulated two results describing the typical behavior and generecity of ${\mathcal M}$-
mean curvature lines for a function ${\mathcal M}= 2{\mathcal H}- {\mathcal M_1}\sqrt{\mathcal K}$.

\begin{proposition}\label{prop:p2b} Let $0$ be a parabolic point and the Monge chart $(x,y)$ as in equation \ref{eq:ez2}.
 Suppose that $\mathcal M $  is  a mean curvature function {\it
1/2-regular}, case b).
\begin{itemize}
\item[i)]
If $b\ne 0$ then the $ {\mathcal M}$-mean  curvature lines are
transversal to the parabolic curve, as  shown in  Figure
\ref{fig:mpr},
 the cuspidal case.

\item[ii)] If $b=0$ and $\sigma= (cd-2d^2+kC)kc m_0^2\ne 0$ then the  $ {\mathcal M}$-mean
curvature lines are shown in the Figure \ref{fig:mpr}. In fact, if
$\sigma >0$ then the $ {\mathcal M}$-mean   curvature lines are
folded saddles. Otherwise,  if $\sigma <0$ then the$ {\mathcal
M}$- mean curvature lines are folded nodes.

\end{itemize}

\end{proposition}
\begin{proof} Similar to that of proposition \ref{prop:p2}.
\end{proof}

\begin{theorem}\label{th:51b}
Assume that the ${\mathcal M}$- mean curvature function $\mathcal M$ is $1/2$-regular as in definition \ref{def:mb2}, case b).

An immersion $\alpha\in {\mathcal I}^{r,s}({\mathbb V^2)}$, $ r\geq   6$, is  $C^6-$local $ {\mathcal M}$- mean curvature structurally stable at a tangential parabolic point $p$ if only if, the condition $\sigma \neq 0$ in  proposition \ref{prop:p2b}  holds.

\end{theorem}

\begin{proof} Similar to that of theorem \ref{th:51a}.
\end{proof}

\section {On $\mathcal M$- Mean Curvature Structural Stability}
\label{sec:6}

In this section the results of sections \ref{sec:3}, \ref{sec:4} and \ref{sec:5}   are put together
 to provide sufficient conditions for {$\mathcal M$-mean curvature stability.

\begin{theorem} \label{th:sta} Let $\mathcal M$ be a mean curvature function which is  positive regular,  1-regular or 1/2-regular.  See definitions \ref{def:bar},  \ref{def:mb1} and \ref{def:mb2}.

Then the set of immersions ${\mathcal G}_i({\mathbb V^2}), i=1,\; 2$ which satisfy
conditions $i$), ... , $v$) below
 are i-$C^s$-$\mathcal M$-mean   curvature structurally stable and ${\mathcal G}_i, i=1,\; 2$ is open in ${\mathcal I}^{r,s}({\mathbb V^2)}, \; r\geq  s\geq 6$.
\begin{itemize}

\item[i)]  The parabolic curve is 1-regular or 1/2-regular: ${\mathcal K}=0$ implies $d{\mathcal K}   \neq 0$
and the tangential singularities are folded saddles,  nodes and foci.

\item[ii)]  The umbilic points are of type $M_i$, $i=1,\; 2,\;3$.

\item[iii)]  The    cycles of $\mathbb H_{\alpha,i}^{\mathcal M}$ are hyperbolic.

\item[iv)]  The    foliations $\mathbb H_{\alpha,i}^{\mathcal M}$  have no separatrix connections. This means that there is no   $\mathcal M$-mean curvature line joining two umbilic or tangential parabolic singularities and being separatrices at both ends. See  propositions \ref{prop:2}, \ref{prop:p1},  \ref{prop:p2}, \ref{prop:p1b} and \ref{prop:p2b}.

\item[v)] The limit set of every leaf of $\mathbb H_{\alpha,i}^{\mathcal M}$ is    a parabolic point, umbilic point or a  $\mathcal M$-mean curvature cycle.
\end{itemize}

\end{theorem}

\begin{proof} The openness of ${\mathcal G}_i({\mathbb V^2})$  follows from the local structure of the  $\mathcal M$-mean curvature lines near the umbilic points of types $M_i$, $i=1,2,3$, near the parabolic points (cusp, saddles, foci and nodes), near the  $\mathcal M$-mean curvature cycles and by the absence of umbilic  $\mathcal M$-mean curvature separatrix connections and the absence of recurrences.
 The  equivalence can be performed by the method of canonical regions and their continuation as was done in
\cite{gs1, gs2} for principal lines,  and in \cite{a2}, for asymptotic lines. \end{proof}

 Notice that   Theorem \ref{th:sta} can be reformulated
so as to give the $\mathcal M$-mean    stability of the configuration rather than
that of the  separate  foliations. To this end it is necessary to consider the folded extended lines, that  is  to consider the line  of one foliation that arrive  at the parabolic set at a given transversal point as continuing through the line  of the other foliation leaving the parabolic set at  this point, in a sort of ``billiard". This gives raise to the extended folded cycles and separatrices that must be preserved by the homeomorphism mapping simultaneously the two foliations.

Therefore the third, fourth and fifth hypotheses above should be modified as follows:
\begin{itemize}
\item[iii')] the extended folded periodic cycles should be hyperbolic,
\item[iv')] the extended folded separatrices should be disjoint,
\item[v')] the limit set of extended lines should be umbilic points, parabolic singularities and extended folded cycles.
\end{itemize}

 The class of immersions which verify the extended five conditions i), ii), iii'), iv'), v') of a compact and oriented manifold $\mathbb V^2$ will be denoted by ${\mathcal G}({\mathbb V^2})$.

This procedure has been adopted by the authors in the case of asymptotic lines by the suspension operation in order to pass from the foliations to the configuration and properly formulate the stability results.  See \cite{a2}.

\begin{remark} In the space   of convex  immersions ${\mathcal I}^{r,s}_c({\mathbb S}^2)$ ( ${\mathcal K}_\alpha>0$),  the sets  ${\mathcal G}({\mathbb S^2})$ and ${\mathcal G}_1({\mathbb S^2})\cap {\mathcal G}_2({\mathbb S^2}) $ coincide. \end{remark}

 The density result involving the five conditions above  is formulated now.

\begin{theorem}\label{th:d} Let $\mathcal M$ be a mean curvature function which is  positive regular,  1-regular or 1/2-regular. See definitions \ref{def:bar},  \ref{def:mb1} and \ref{def:mb2}.

Then the sets ${\mathcal G}_i, \; i=1,\; 2$,  are dense in ${\mathcal I}^{r,2}({\mathbb V^2)},\; r\geq  6$.

In the space  ${\mathcal I}^{r,2}_c({\mathbb S}^2)$ the set ${\mathcal G}({\mathbb S^2})$ is dense.
\end{theorem}

The main ingredients for the proof of this theorem are the Lifting
and   Stabilization Lemmas, essential for the achievement of
condition five, are developed in next section.

\section{ Density of  $\mathcal M$-Mean Curvature Structurally  Stable Immersions}\label{sec:7}

In this section will be proved  an approximation theorem for the class of immersions or surfaces
having structurally stable  {\it
$\mathcal M$- mean   curvature configuration}.

The proof of Theorem \ref{th:d}  follows from the elimination of $\mathcal M$- mean   curvature recurrences and
 the stabilization of the $\mathcal M$- mean curvature  separatrices.
 The steps are basically   those followed by C. Gutierrez and J. Sotomayor
in the case of principal curvature lines, see \cite{gsln, gs2}.
The main ideas goes back to M. Peixoto \cite{mp} and C. Pugh \cite{pu} to solve the similar problem of
elimination of recurrences
for vector fields on surfaces. See also the book by W. de Melo and J. Palis \cite{pm}.

In what follows will be established
 the main
 technical lemmas necessary to obtain the Lifting Lemma, essential to control  the  effect on $\mathcal M$-mean   curvature lines under
 suitable deformations of the immersion.
 There is no lost of
generality to assume
that the immersion is $C^\infty$ or $C^\omega$ in the proof of the density theorem.

In what follows a chart whose coordinates lines are  $\mathcal M$-mean   curvature lines will be called  a $\mathcal M$-{\it  mean
curvature chart} for the immersion $\alpha$.

\begin{lemma} \label{lem:71}
Let $\alpha:\mathbb V^2\to \mathbb R^3$ be an immersion of class $C^\infty$  and $(u,v):(U,D)\to (V,I\times I)$
be a positive $\mathcal M$-mean
curvature chart on $\mathbb V^2$, where $I=[-1,1]$. Suppose that, for $\ep$ small,
$\beta=\alpha_\ep=\alpha+\ep\varphi
N $ is an immersion and $\varphi$ be a smooth function on $U$ which satisfies:
$\varphi(-1,v) =\varphi(1,v)=\varphi_u(-1,v)=\varphi_u(1,v)=\varphi_{uu}(-1,v)=\varphi_{uu}(1,v)=0$.
Then the $\mathcal M$-mean curvature line of $\alpha_\ep$ on $D$
 which passes through $q$ in $\{u=-1\}\cap \{-1<v<1\}$ meets the segment of abscissa $\{u=1\}$ at a point  whose $v-$coordinate $v_\ep$ has a derivative with respect to $\ep$ given by:

\begin{equation}\label{eq:71}
\aligned \frac{d}{d\ep}(v_\ep)|_{\ep=0}
 &=\int_{-1}^1  \frac{E[ 2{\mathcal M}{\mathcal M}_K + {\mathcal M}_H]}{4\tau_g \sqrt{EG-F^2} } \varphi_{vv}du+\int_{-1}^1 A_1(u)\varphi_vdu\\
&+\int_{-1}^1A_2(u)\varphi du\endaligned
\end{equation}
where $A_1$ and $A_2$ are functions of the coefficients of the first and second fundamental form of $\alpha$.
\end{lemma}

\begin{proof} Suppose that for $\epsilon$ small, $$\beta(u,v,\epsilon)=\alpha_\ep(u,v)=\alpha(u,v)+
\epsilon \varphi(u,v)N(u,v)$$
\noindent  is an immersion.

The $v$-coordinate, $\; v=v(u,q,\ep)$, of the point where the line of $\mathcal M$-mean
 curvature through the point $q$
 in $\{u=-1\}\cap \{-1<v<1\}$  meets the curve with abscissa $\{u\}$,
 satisfies the following Cauchy Problem with parameter $\epsilon$.

\begin{equation}\label{eq:gm}
 (e-{\mathcal M}E)+ 2(f-{\mathcal M}F)\frac{dv}{du} +
 (g-{\mathcal M}G)(\frac{dv}{du})^2 =0,
\quad v(-1,\epsilon)=q
\end{equation}

Since  $(u,v)$ is a  $\mathcal M$-mean curvature chart,  it results that

\begin{equation}\label{eq:mcc}
\aligned
 \frac{dv}{du}(u,q,0)=&0,\\
(e-{\mathcal M}E)(u,v,0)=& (g-{\mathcal M}G)(u,v,0)=0,\;\;\\
 (f-{\mathcal M}F)(u,v,0)=& \tau_g \sqrt{EG-F^2}\\
=& \sqrt{(k_2-{\mathcal M})({\mathcal M}-k_1)}\sqrt{EG-F^2}\ne 0
\endaligned
\end{equation}

Differentiating the equation (\ref{eq:gm}) with respect to $\epsilon$,  evaluated
 on $(u,v(q),\ep)$, making $\ep=0$ and using (\ref{eq:mcc}) it follows that
$$\frac{dv_\ep}{d\ep}=\frac{\partial v_\ep}{\partial \ep}(u,q,\ep)|_{\ep=0}$$
\noindent  satisfies
the following Cauchy Problem:

\begin{equation} \label{eq:74} \aligned \frac{d }{du}&(\frac{dv_\ep}{d\ep})
 =-\frac{[  e_\ep-  {\mathcal M}_\ep E-  {\mathcal M}E_\ep ]}{2(f- {\mathcal M}F)}
 (u,v(q),0)\\
 \frac{dv_\ep}{d\ep}&(-1,q,0)=0\endaligned \end{equation}

The structure  equations for the immersion  $\alpha$ are given by:

\begin{equation}
\label{eq:75} \aligned
 N_u &=\frac{fF-eG}{EG-F^2}\alpha_u+\frac{eF-fE}{EG-F^2}\alpha_v\\
N_v &=\frac{gF-fG}{EG-F^2}\alpha_u+\frac{fF-gE}{EG-F^2}\alpha_v\\
\alpha_{uu} &= \Gamma_{11}^1\alpha_u+\Gamma_{11}^2\alpha_v+eN\\
\alpha_{uv} &= \Gamma_{12}^1\alpha_u+\Gamma_{12}^2\alpha_v+fN\\
\alpha_{vv} &= \Gamma_{22}^1\alpha_u+\Gamma_{22}^2\alpha_v+gN\endaligned\end{equation}

The functions $\Gamma_{ij}^k$ are the Christoffel symbols whose expressions in terms of $E$, $F$  and $G$ in a   chart are $(u,v)$ are given by:
\begin{equation}\label{eq:ch}\aligned \Gamma_{11}^1&=\frac{GE_u-2FF_u+FE_v}{2(EG-F^2)}, \;\;\;\;
\Gamma_{11}^2=\frac{2EF_u-EE_v-FE_u}{2(EG-F^2)},\\
\Gamma_{12}^1&=\frac{GE_v-FG_u}{2(EG-F^2)},\hskip 2cm
\Gamma_{12}^2=\frac{EG_u-FE_v}{2(EG-F^2)},\\
\Gamma_{22}^1&=\frac{2GF_v-GG_u-FG_v}{2(EG-F^2)},\;\;\;\;
\Gamma_{22}^2=\frac{EG_v-2FF_v+FG_u}{2(EG-F^2)}.
\endaligned
\end{equation}

By direct calculation, it is obtained that

\begin{equation}\label{eq:77}
\aligned \beta_{u}&=(1 +\ep\varphi\frac{fF-eG}{EG-F^2})\alpha_u+\ep\varphi\frac{eF-fE}{EG-F^2}   \alpha_v+
\ep\varphi_u N\\
\beta_v &=  \ep\varphi \frac{gF-fG}{EG-F^2}   \alpha_u+(1+\ep\varphi\frac{fF-gE}{EG-F^2} )\alpha_v+\ep \varphi_v N\endaligned \end{equation}

\begin{equation}\label{eq:78}
\aligned
\beta_{uu}  =&[ \Gamma_{11}^1+\ep\varphi \big(\Gamma_{11}^1\frac{fF-eG}{EG-F^2} +\Gamma_{12}^1\frac{eF-fE}{EG-F^2}+(\frac{fF-eG}{EG-F^2})_u \big)\\
+&2\ep\varphi_u \frac{fF-eG}{EG-F^2}  ]\alpha_u\\
 +&[\Gamma_{11}^2+\ep\varphi \big(\Gamma_{11}^2\frac{fF-eG}{EG-F^2} +\Gamma_{12}^2\frac{eF-fE}{EG-F^2}+(\frac{eF-fE}{EG-F^2})_u \big)\\
+& 2\ep\varphi_u \frac{eF-fE}{EG-F^2}  ]\alpha_v\\
 +&[e+\ep\varphi\frac{2efF-e^2G-f^2E}{EG-F^2}+\ep\varphi_{uu}]N\endaligned
\end{equation}

\begin{equation}\label{eq:79}
\aligned \beta_{uv} =&
[ \Gamma_{12}^1+\ep\varphi \big(\Gamma_{12}^1\frac{fF-eG}{EG-F^2} +\Gamma_{22}^2\frac{eF-fE}{EG-F^2}+(\frac{fF-eG}{EG-F^2})_v \big)\\
+& \ep\varphi_u \frac{gF-fG}{EG-F^2}  + \ep\varphi_v \frac{fF-eG}{EG-F^2} ]\alpha_u\\
 +&[\Gamma_{12}^2+\ep\varphi \big(\Gamma_{12}^2\frac{fF-eG}{EG-F^2} +\Gamma_{12}^2\frac{eF-fE}{EG-F^2}+(\frac{eF-fE}{EG-F^2})_v \big)\\
+&  \ep\varphi_u \frac{fF-gE}{EG-F^2} +  \ep\varphi_v \frac{eF-fE}{EG-F^2} ]\alpha_v\\
 +&[f+\ep\varphi(f\frac{ fF-e G }{EG-F^2}+g\frac{eF-fE}{EG-F^2}+\ep\varphi_{uv}]N\endaligned
 \end{equation}

\begin{equation}\label{eq:710}
 \aligned \beta_{vv} =&
[ \Gamma_{22}^1+\ep\varphi \big(\Gamma_{22}^1\frac{fF-gE}{EG-F^2} +\Gamma_{12}^1\frac{gF-fG}{EG-F^2}+(\frac{gF-fG}{EG-F^2})_u \big)\\
+&2\ep\varphi_v \frac{gF-fG}{EG-F^2}  ]\alpha_u\\
 +&[\Gamma_{22}^2+\ep\varphi \big(\Gamma_{22}^2\frac{fF-gE}{EG-F^2} +\Gamma_{12}^2\frac{gF-fG}{EG-F^2}+
(\frac{fF-gE}{EG-F^2})_v \big)\\
+& 2\ep\varphi_v \frac{fF-gE}{EG-F^2}  ]\alpha_v\\
 +&[g+\ep\varphi(g\frac{ fF-gE   }{EG-F^2}+f\frac{gF-fG}{EG-F^2})+\ep\varphi_{vv}]N\endaligned
\end{equation}
Also,
\begin{equation}\label{eq:711}
\frac{\partial}{\partial\ep}(|\beta_u\wedge \beta_v|)|_{\ep=0}
=-2\varphi{\mathcal H}\sqrt{EG-F^2}
 \end{equation}

Therefore, using the equations (\ref{eq:77}) - (\ref{eq:710}) the following  is obtained.
\begin{equation}\label{eq:712}
\aligned E_\ep &=-2\varphi e,\qquad\qquad
F_\ep =-2\varphi f,\qquad\qquad
G_\ep =-2\varphi g\\
e_\ep &=\varphi_{uu}+\varphi [\frac{2efF-e^2G-f^2F}{EG-F^2}] -   \varphi_u\Gamma_{11}^1-\varphi_v\Gamma_{11}^2
 \\
f_\ep &=\varphi_{uv}+\varphi [f\frac{fF-eG}{EG-F^2}+g\frac{eF-fE}{EG-F^2}] -   \varphi_u\Gamma_{12}^1-\varphi_v\Gamma_{12}^2 \\
g_\ep &=\varphi_{vv}+\varphi [g\frac{fF-gE}{EG-F^2}+f\frac{gF-fG}{EG-F^2}] -   \varphi_u\Gamma_{22}^1-\varphi_v\Gamma_{22}^2.\endaligned \end{equation}
Let
\begin{equation}\label{eq:76} {\mathcal H}=\frac{Eg+eG-2fF}{2(EG-F^2)}, \qquad  {\mathcal K}=\frac{eg-f^2}{ EG-F^2 }. \end{equation}

Then, using equations  (\ref{eq:712}) and (\ref{eq:76}), it follows that

\begin{equation}\label{eq:ke}
\aligned {\mathcal K}_\ep (EG-F^2) =& \varphi K_1 +\varphi_u (2\Gamma_{12}^1f-g\Gamma_{11}^1-e\Gamma_{22}^1)\\
+& \varphi_v(2\Gamma_{12}^2 f-g\Gamma_{11}^2-e\Gamma_{22}^1)+g\varphi_{uu}-2f\varphi_{uv}+e\varphi_{vv}\\
 2{\mathcal H}_\ep (EG-F^2) =& \varphi H_1 -\varphi_u  \Gamma_{22}^1 -
  \varphi_v \Gamma_{22}^1  +\varphi_{vv}\endaligned
\end{equation}

where,

$$K_1=K_1(e,f,g,E,F,G)(u,v,0), \; \; H_1=H_1(e,f,g,E,F,G)(u,v,0).$$

Now as ${\mathcal M}=m({\mathcal H},{\mathcal K})$ it follows,  from equations (\ref{eq:mcc}) and (\ref{eq:ke}), that

$$\aligned \frac{d}{d\ep}{\mathcal M}_\ep =&{\mathcal M}_H {\mathcal H}_\ep + {\mathcal M}_K {\mathcal K}_\ep \\
=& \frac{\varphi_{vv}}{EG-F^2}( \frac{{E\mathcal M}_H}{2} + E{\mathcal M}{\mathcal M}_K)+ (.) \varphi +(.)\varphi_u+(.) \varphi_{v}+(.)\varphi_{uv}   \endaligned $$

Therefore,

\begin{equation}\label{eq:713}
\aligned \frac{d}{d\ep}\big(\frac{ e-{\mathcal M}E}{2(f-{\mathcal M}F)}\big)|_{\ep=0}=& \frac{E(  {{\mathcal M}_H} + 2{\mathcal M} {\mathcal M}_K) }{4(EG-F^2)(f-{\mathcal M}F)} \varphi_{vv}\\
+& (.) \varphi +(.)\varphi_u+(.) \varphi_{v}+(.)\varphi_{uv} +(.)\varphi_{uu}
\endaligned
\end{equation}
Here $(.)$ denote  functions involving   $\mathcal M(\mathcal H, \mathcal K)$ and the coefficients of the first and second fundamental forms of $\alpha$.

Using (\ref{eq:713}) when integrating the variational equation (\ref{eq:74}) and performing the partial integration
 with  boundary conditions on the function $\varphi$,  the expression for
$(\frac{dv_\ep}{d\ep})|_{\ep=0}$ is achieved, as stated in (\ref{eq:71}).
\end{proof}

\begin{lemma}\label{lem:72} Let $\alpha:\mathbb V^2\to \mathbb R^3$ be an immersion of class $C^\infty$ and
$(u,v):(U,D)\to (V,I\times I)$ be a positive ${\mathcal M}$-mean curvature chart on $\mathbb V^2$, where $I=[-1,1]$.
  Suppose that ${{\mathcal M}_H} + 2{\mathcal M} {\mathcal M}_K>0 $, i.e, $\mathcal M$ is positive regular. Then there exists a smooth function $\varphi:\mathbb V^2\to [0,1]$ whose support is contained in $D$ such that, if
 $\ep$ is small enough then, for every $\ep$ in $[-r,r]$, $\beta=\alpha+\ep\varphi N $ is an immersion and the ${\mathcal M}$-mean
 curvature line for $\beta $ on  $D$ which passes through $q$ in $\{u=-1\}\cap \{-1<v<1\}$ meets the segment
$\{u=1\}\times \{-1<v<1\}$ at a point $v_\ep(q)$ so that the map $\ep\to v_\ep(q)$ is strictly increasing.
\end{lemma}

 \begin{proof} Let $\rho$ be a real smooth function with values in $[0,1]$, identically equal to $1$ on a neighborhood
of $0$ and with support contained in $I$.

 Let  $\varphi=\varphi(u,v)=b\frac{v^2}2\rho(u)\rho(v)$ and take $r>0$ small so that for any $\ep$ in
$[-r,r]$,
 $\beta=\alpha_\ep =\alpha+\ep\varphi N$ is a smooth  immersion. Let $v(u)$, $u\in I$, be the $v-$coordinate of the ${\mathcal M}$-mean
 curvature lines of $\alpha_\ep=\alpha+\ep\varphi N$ such that $v_\ep(q)=q.$ As $\varphi(u,0)=\varphi_v(u,0)=0$ and
$\varphi_{vv}(u,0)=b\rho(u)$ by lemma \ref{lem:71} applied to the family of immersions $\alpha_\ep$ it follows that

 $$\frac{\partial v}{\partial \ep}(u,\ep)|_{(0,0)}= \int_{-1}^1
\frac{E (  {{\mathcal M}_H} + 2{\mathcal M} {\mathcal M}_K) }{4\tau_g \sqrt{EG-F^2}}
\rho(u)du
= c>0.$$

 \noindent This implies that the map $\ep\to v_\ep(q)$ is strictly increasing. This proves the lemma.\qquad \qed
\end{proof}

\begin{lemma}\label{lem:73} Let $\alpha:\mathbb V^2\to\mathbb R^3$ be an immersion  and $(u,v):(U,D)\to (V,I\times I)$
be a positive ${\mathcal M}$-mean curvature chart for $\alpha$ on $\mathbb V^2$, where $I=[-1,1]$. Assume also that $\mathcal M$ is positive regular.
 Then given any $\eta >0$, there are  numbers $d,c\in (0,\frac 1{12})$ such that for
every $r\in (0,d]$ and $q$ in $\{u=-1\}\cap \{-\frac 12<v<\frac 12\}$,
 there exists a smooth function $\varphi:\mathbb V^2\to [0,1]$ whose support is
contained in $D_r=v^{-1}(v(q)+rI)$ and $||\varphi||_{2,V}$, the $C^2-$norm of $\varphi$ on
 $V$, in the $(u,v)-$coordinate chart, is less than $\eta$.

Furthermore, for every   $\ep\in I$, $\alpha_\ep=\alpha+\ep\varphi N $ is an
immersion and the ${\mathcal M}$-mean curvature line for $\alpha_\ep$ on  $D$ which passes through $q$
 in $\{u=-1\}\cap \{-1<v<1\}$ meets the segment  $\{u=1\}\cap \{-1<v<1\}$ at a point $v_\ep(q)$ so
that the map $\ep\to v_\ep(q)$ is strictly increasing and its image
contains the interval $[v(q)-2c\ep,v(q)+2c\ep]$.
\end{lemma}

\begin{proof} Let $\rho$ be a real smooth function with values in $[0,1]$, identically equal to $1$
 on $\frac 56 I$ and with support contained in $\frac 67 I$. Let also $\eta >0 $ be given.
There are real numbers $c>0$ and $b$ such that for all $(u_0,v_0)$ in $I\times I$ it  follows that,

\begin{equation}\label{eq:714a}  6|b|(||\rho||_2)^2<\eta \end{equation}

 \begin{equation}\label{eq:714b} \int_{-1}^0 b
\frac{E  (  {{\mathcal M}_H} + 2{\mathcal M} {\mathcal M}_K)(u,v_0) }{4\tau_g \sqrt{EG-F^2}}\rho(u)du
<\frac 14\end{equation}

\begin{equation}\label{eq:714c}
\int_{-1}^1 b
\frac{E (  {{\mathcal M}_H} + 2{\mathcal M}
 {\mathcal M}_K)(u,v_0)  }{4\tau_g \sqrt{EG-F^2}}
\rho(u)du >
3c \end{equation}

Let $\psi$ be a smooth real function on $U\times I\times I$, defined by

\centerline{$\psi(u,v;v_0,\ep)=b\ep\frac{(v-v_0)^2}2 \rho(u).$}

It will be proved that if $d=d(\eta)\in (0,\frac 1{12})$ is small enough, then for
 every $\ep\in(0,d)$ and $q$  with $u(q)=-1$ and $v(q)=v_0$ in $\frac 12 I$,
 the smooth function $\varphi(.)=\varphi(.;v(q),\ep)$ defined on $\mathbb V^2$ by

$$\varphi(u,v;v_0,\ep)=\psi(u,v;v_0,\ep)\rho(\frac{v-v_0}{|\ep|})$$
\noindent  whose support is contained in $D_r$, satisfies the conditions required by the lemma.

In fact, suppose $d>0$ is so small that for any $(v_0,\ep)\in I\times dI$,
 $\alpha_{v_0,\ep}=\alpha+\psi(.;v_0,\ep)N$ is an immersion.

Let $v(u;v_0,\ep)$, $u\in I$ and $v\in \frac 45 I$, be the $v-$coordinate of
 a ${\mathcal M}$-mean curvature line of $\alpha_\ep$ through the point $q$, with $u(q)=-1$
 and $v(q)=v_0$. As $\psi(u,v_0;v_0,\ep)=0$, using (\ref{eq:714a}), (\ref{eq:714b}) and (\ref{eq:714c}), it  follows
 from lemmas
\ref{lem:71}  and \ref{lem:72},
applied to the family of immersions $\alpha_{v_0,\ep}$, depending on the parameter $\ep$, that for all $(u,v_0)$ in $I\times \frac 34 I$,

\centerline{$\frac{\partial v}{\partial \ep}(1,v_0,0)>2c \qquad \qquad\text{and}\qquad\qquad
\frac{\partial v}{\partial \ep}(1,v_0,0)<\frac 13.$}

\noindent Hence, as $I $ is compact and $ \frac{\partial v}{\partial \ep}(1,v_0,\ep)$ depends continuously on $(v_0,\ep)$, taking $ d>0 $ small enough, it holds that for all $(u;v_0,\ep)$ in $I\times\frac 12 I\times dI$,

\centerline{$\frac{\partial v}{\partial \ep}(1,v_0,0)>c \qquad\qquad \text{and}\qquad\qquad
\frac{\partial v}{\partial \ep}(1,v_0,0)<\frac 12.$}

Therefore from the Mean Value Theorem, for all $(u;v_0,\ep)$ in $I\times \frac 12 I\times dI$

\begin{equation}\label{eq:715abc}
\aligned     v&(1;v_0,\ep)\geq v_0+c\ep,\;\;\text{ if}\;\; \ep\leq 0    \\
  v&(1;v_0,\ep)\leq v_0+c\ep, \;\;\text{if}\;\; \ep <0\;\;\; \text{and} \\
 |&v(u;v_0,\ep)-v_0|\leq \frac{|\ep |}2    \endaligned      \end{equation}

Without lost of generality, $\eta>0$
can be assumed to be so small that for $\ep\in I$, $\alpha_\ep=\alpha+\ep\varphi N$ is an immersion.
 Notice that for all $v$ in $v(q)+\frac 56 \ep I$, $\rho(\frac{v-v_0}{|\ep|})=1$.
 Therefore for all $\ep $ in $I$ and all $(u,v)$ in $I\times (v(q)+\frac 56\ep I)$,
 $\varphi(u,v;v_0,\ep)=\psi(u,v;v_0,\ep)$. From this and $(7.15c)$ follows that the ${\mathcal M}$-mean   curvature
lines of $(\alpha+\psi(.,v_0,\ep)N)|_{D^\prime}$ and those of the
$\alpha+\varphi(.,v_0,\ep)N|_{D^\prime}$ coincide, where $D^\prime=v^{-1}[v(q)+\frac 56 \ep I]$.
Hence, the assertion that the map $v_\ep(q)$ is strictly increasing
 and its image contains the interval $[v(q)-2c\ep,v(q)+2c\ep]$ follows from equation  \ref{eq:715abc}.
\end{proof}

\begin{lemma}\rm{[{\bf Lifting Lemma}]}\label{lem:74}
Let $\alpha:\mathbb V^2\to \mathbb R^3$ be an immersion of class $C^\infty$  with a minimal  ${\mathcal M}$-mean curvature line
$\tilde\gamma$, oriented
from a starting  point $q$, whose $\omega -$limit set contains a nontrivial minimal  recurrent ${\mathcal M}$- mean
curvature line  $\gamma$. Assume also that $\mathcal M$ is positive regular. Then given any $\eta>0$, $p\in\gamma$ and any ${\mathcal M}$-mean curvature chart
$(u,v):(U,D,p)\to  (V, I\times I,0)$ where $I=[-1,1]$, there is a ${\mathcal M}$-mean  curvature chart
$(s,t):D^\prime \to I\times I$ for $\alpha$ on $\mathbb V^2$ and a smooth function $\varphi:\mathbb V^2\to [0,1]$
such that:
\item{i)} the support of $\varphi$ is contained in $D\cap D^\prime$ and $||\varphi||_{2,V}$,
 the $C^2-$norm of $\varphi$ on $V$,   is less than $\eta$.

\item{ii)} There are arcs of minimal  ${\mathcal M}$-mean   curvature lines $[b,a]\subset [q,a]\subset \tilde\gamma$ such that $a$, $b$ are in the arc $\{s=-1\}$ and $[b,a]\cap D^\prime=[q,a]\cap D^\prime \subset [a,a^\prime]\cup [b,b^\prime]$, where $a^\prime$ and $b^\prime$ are the points on $\{u=1\}$, defined by $t(a^\prime)=t(a)$ and $t(b^\prime)=t(b)$.

Moreover the minimal  ${\mathcal M}$-mean   curvature lines for $\alpha_\ep$ on
 $D^\prime$ which passes through $a$ ( resp. $b$) meets the segment $\{u=1\}$ at a point $v_\ep(a)$ ( resp. $v_\ep(b)$) in such way that for some values of $\ep\in [0,1]$, it coincides with $a^\prime$ and $b^\prime$. See Figure  \ref{fig:71}.
\end{lemma}
\begin{proof}
See     \cite{gsln, gs2}.\qquad
 \end{proof}

\begin{proposition}\label{prop:71} Let $\alpha:\mathbb V^2\to \mathbb R^3$ be an immersion of class $C^\infty$ with
umbilic set ${\mathcal U}_{\alpha}\ne\emptyset$ and having a nontrivial minimal  recurrent
 ${\mathcal M}$-mean curvature line $\gamma$ and let $A$ be a subset of $\mathbb V^2$ formed by finitely many
 minimal  ${\mathcal M}$-mean curvature lines that are either minimal ${\mathcal M}$-mean curvature  separatrix connections
 or minimal ${\mathcal M}$-mean   curvature cycles. Assume also that $\mathcal M$ is positive regular. Then there is a point $p\in \gamma\setminus A$ such
 that given any chart $(u,v):(U,p)\to (V,0)$ on a neighborhood $U$ of $p$, where $U $ is
 disjoint of $A$, there is a sequence of smooth functions $\varphi_n$ on $\mathbb V^2$,
 whose support is contained in $U$ such that $||\varphi_n||_{2,V}$,
 the $C^2-$norm of $\varphi_n$, in the coordinate chart $(u,v)$,  tends to $0$ and
such that the immersions $\alpha_n=\alpha+\varphi_n N$ satisfy the following alternatives:
\item{i)} $\alpha_n$ has a ${\mathcal M}$-mean   curvature cycle $\gamma_n$ not completely contained
 in $\mathbb V^2\setminus U$. Moreover if there is a minimal  ${\mathcal M}$-mean   curvature cycle of
 $\alpha$ (i.e. disjoint of $U$) which together with $\gamma_n$ bound a cylinder in
$\mathbb V^2$ then this cylinder contains an umbilic point of $\alpha$.
\item{ii)} $\alpha_n$ has at least one minimal  ${\mathcal M}$-mean   curvature separatrix connection more than the
immersion $\alpha$ does.
\end{proposition}

\begin{proof} See    \cite{gsln, gs2}.
\end{proof}

\begin{proposition}\label{prop:72} Let $\alpha:\mathbb V^2\to \mathbb R^3$ be an immersion of class $C^\infty$
 and let $A$ be a subset of $\mathbb V^2$ formed by finitely many minimal   ${\mathcal M}$-mean   curvature lines
that are either minimal  ${\mathcal M}$-mean   curvature separatrix connections or minimal   ${\mathcal M}$-mean
curvature cycles. Assume also that $\mathcal M$ is positive regular.

 Then there is a sequence of immersions $\alpha_n=\alpha+\varphi_nN$, $C^2$-converging to $\alpha$, such that
 the support of $\varphi_n$ is disjoint from $\bar A=A\cup {\mathcal U}_{\alpha}$ and $\alpha_n$ has
 no non trivial minimal  recurrent ${\mathcal M}$- mean   curvature lines.\end{proposition}

\begin{proof} See     \cite{gsln, gs2}.
\end{proof}

Let $\alpha:\mathbb V^2\to\mathbb R^3 $ be an immersion whose
umbilic points are of type $M_i$, $\; i=1,2,3$ and all tangential
parabolic singularities are folded saddles, nodes and foci. A
minimal  ${\mathcal M}$- mean   curvature separatrix $\Gamma$ of
an umbilic or parabolic point $p$ is said to be {\it stabilized}
provided:
\begin{itemize}
\item[i)] it is not a minimal  ${\mathcal M}$- mean   curvature separatrix connection for $\alpha$;
\item[ii)] its limit sets are umbilic points, parabolic point or attracting or repelling minimal   ${\mathcal M}$- mean   curvature
cycles, and
\item[iii)] $\alpha$ is in the $C^6$-interior of the set of immersions that satisfy i) and ii); i.e.,
for any sequence of
immersions $\alpha_n$, $C^6$-converging to $\alpha_n$, the sequence of separatrices $\Gamma_n$, of an  umbilic or a tangential parabolic point $p_n$,
converging to the separatrix $\Gamma$ of $p$, verify i) and ii) for $\alpha_n$.
\end{itemize}

\begin{lemma}\label{lem:75} {\rm\bf [Stabilization Lemma]} Any immersion  $\alpha:\mathbb V^2\to \mathbb R^3$ of class
$C^\infty$  is the
$C^2$-limit of a sequence of immersions whose umbilic points are all of type $M_i$, $i=1,2,3$, section \ref{sec:3}, all tangential parabolic singularities are folded saddles, nodes and foci, and furthermore:
\begin{itemize}
\item[i)] their minimal ${\mathcal M}$-mean   curvature   separatrices are all stabilized;
\item[ii)] the $\omega$-limit set of any oriented minimal ${\mathcal M}$-mean curvature line is either an umbilic or a cuspidal parabolic point   or a minimal ${\mathcal M}$-mean   curvature cycle, and
\item[iii)] for any $s\geq 6$, $\alpha$ is in the $C^s$-interior of the set of immersions
 satisfying i) and ii).
\end{itemize}
\end{lemma}

\begin{proof} See       \cite{gsln, gs2}.\qquad
\end{proof}

\begin{remark}\label{rem:71} In all lemmas and propositions above the same
 conclusions  hold for the maximal ${\mathcal M}$-mean   curvature lines provided the corresponding hypotheses are
made also in this case.
\end{remark}

\subsection{  Proof  of the Density Theorem \ref{th:d}}\hskip 8cm
  \newline

\noindent { \bf Part 1: Elimination of nontrivial recurrences}

By proposition \ref{prop:72} the recurrent lines can always be destroyed
by a finite sequence of small local $C^2-$perturbations of the
immersion $\alpha$.
Each perturbation creates either a new ${\mathcal M}$-mean   curvature
cycle or a new $\mathcal M$- mean   curvature separatrix connection.

Initially will be considered
  the elimination of the minimal recurrent $\mathcal M$- mean curvature
lines.

The key points involved in the argument will be
 given below.

Let $\gamma $ be a non trivial minimal recurrent $\mathcal M$- mean curvature line.

Assume first that  $\gamma $ is {\it orientable}, i.e,
it  is possible to give an orientation
in
$\gamma$ such
that on a $\mathcal M$- mean   curvature chart it is induced by an   orientation
defined locally  on  the $\mathcal M$- mean curvature line field by the chart.
The recurrent lines on vector fields and those of $\mathcal M$- mean   curvature
foliations on the Torus, in section 5,
are of this type.
In this case  there is a
piecewise smooth simple closed curve of the form $[\b b,\b
a]\cup[\b b;\b a]$, with $[a,b]\subset \gamma $
  and $\b a$  near $\b b,$ that can be
slightly perturbed to
obtain a minimal $\mathcal M$- mean curvature  cycle for the approximating immersion.
Here, and in what follows, $[\b b;\b a]$ means an arc of a maximal  $\mathcal M$- mean curvature line and
 $[\b b,\b a] $ is an arc of a minimal $\mathcal M$- mean curvature line. The arrangement of these points are illustrated in Figure \ref{fig:71}.a.

\begin{figure}[htbp]
 \begin{center}
 \includegraphics[angle=0, width=9cm]{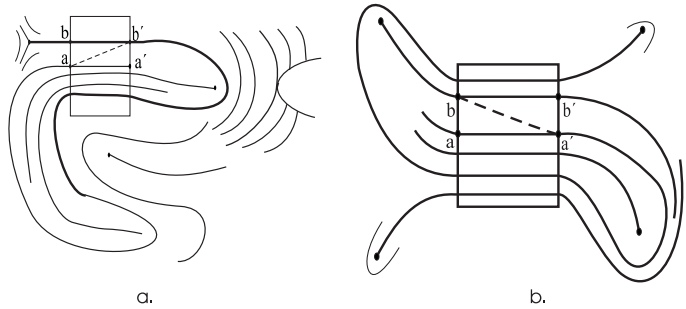}
 \caption{     \label{fig:71}     Recurrences of $\mathcal M$- mean   curvature lines }
    \end{center}
  \end{figure}

When the recurrence is {\it  oscillatory}  ( i.e.  non--orientable),
 then  there is no such
simple closed curve available. In  this case  there are minimal $\mathcal M$- mean   curvature
separatrices accumulating on $\b p. $
These separatrices can be
connected by means of a small perturbation of the immersion. The $\mathcal M$- mean   curvature lines, for ${\mathcal M}={\mathcal H}$, ${\mathcal M}=\sqrt{\mathcal K}$ and ${\mathcal M}={\mathcal K}/{\mathcal H}$,  on the
ellipsoid,    presents this type of
recurrence, see \cite{m}, \cite{g} and \cite{h}.
This situation
is illustrated in Figure \ref{fig:71}. b.

The possibility of finding perturbations as those described
above is established in the Lifting Lemma  \ref{lem:74} and Proposition \ref{prop:71}.

This is done as follows.

Consider a non trivial recurrent minimal separatrix
$\gamma'$ of an umbilic or parabolic point $ q$.
Take $ p\in\gamma$ and a $\mathcal M$- mean curvature chart
$(u,v):(D, p)\to(I\times I,0)$.

By lemma \ref{lem:71}, arbitrarily close to $\{v=0\}$, two points $ \b a$,\; $\b b$  in
$\{u=-1\}\cap\gamma'$ can be selected such that $v(\b b)-v(\b a)=2r>0$
and $(\b b,\b a)$ has
the following spacing property relatively to the  maximal $\mathcal M$- mean
curvature arc $[\b a;\b b]:$

$(\b b,\b a) $ is disjoint from $\{v(\b a)-(3/2)r <v <v(\b
b)+(3/2)r\}\cap\{u=-1\}.$

It results  from this that a local version
of the  lifting argument  (Lemma \ref{lem:72}) can be applied to obtain,
by
means of an $\epsilon-$ small $C^2$-perturbation supported on
$\{v(\b a)-(3/2)r<v<v(\b b)+(3/2)r\}$, the following.
Given $\eta >0,$ there is a constant
$c=c(\eta,(u,v))>0,$ which does
not depend on how close $\b a$  and $\b b$  are, such that, for
every $\b x\in[\b a;\b b],$ the $\mathcal M$- mean curvature  line through the point $\b x$
can reach any
point of the segment $\{u=1\}\cap\{v(\b x)-2cr< v <v(\b x)+2cr\}$.

Consider first the  assumption that $c=1$.

If $[\b b,\b a]
\setminus\{v(\b a)-\frac{3r}{2} <v<v(\b b)+\frac{3r}{2} \}=(\b b',\b a)$,
with $\b b'\in\{u=1\}$,
then, via a perturbation, a minimal $\mathcal M$- mean   curvature cycle can be obtained.
This is illustrated in Figure   \ref{fig:71}.a.

If, however, $[\b b,\b a]
\setminus\{v(\b a)-\frac{3r}{2} <v<v(\b b)+\frac{3r}{2} \}=(\b
b',\b a')$ with
$\b b',\b a'\in\{u=1\}$, it seems to be difficult to approximate $[\b
b',\b a']$ by
a minimal $\mathcal M$- mean   curvature cycle.
Nevertheless, by moving $\b a'$ towards $\b b'$, one
can generate a continuous family of minimal $\mathcal M$- mean   curvature arcs with
endpoints in $\{u=1\}$.
In this process the resulting endpoints
become close to each other but
cannot coincide.
Using this it is
proved that the limit set of this family of minimal $\mathcal M$- mean   curvature arcs must
contain an umbilic or parabolic point $\b q'$ and an arc $[\b p',\b q']$, of
a minimal $\mathcal M$- mean   curvature
 separatrix of $\b q'$, intersecting

\centerline{$
\{v(\b a)-\frac{3r}{2}  < v < v(\b b)+\frac{3r}{2} \}\cap\{u=1\}
$}

\noindent at a point of $(\b b';\b a')$.
In this situation, via a perturbation, a
minimal $\mathcal M$- mean   curvature separatrix connection between $\b q$  and $\b q'$ can be
produced.
This is illustrated in Figure    ref{fig:71}.b.

In general, $c>0$ is much smaller than 1 and
the analysis is done by showing that lemma \ref{lem:72} can be used a number
$n$ of times, where $n$  is of the order of $1/c$, to finally obtain
enough  lifting as to make possible the application of the
arguments above.
The n intervals $[\b a_i;\b b_i]$ which play the same
role as that performed by $[\b a;\b b]$ and on which Lemma \ref{lem:72} is to be
used, are described below.

If $\b a$  and $\b b$  are close enough to each other and $(\b a,\b
b)$ is long
enough, it is proved that there is a family $\{[\b a_t;\b b_t];
t\in[1,n]\}$
of pairwise disjoint maximal $\mathcal M$- mean curvature  arcs such that:
\begin{itemize}
\item[i)] $[\b a_1;\b b_1]=[\b a;\b b]$,

\item[ii)] the curves $\b a_t\in\tilde\gamma,
\b b_t\in\tilde\gamma$ are regular,

\item[iii)]
for all $i\in\{1,2,...,n\}$, $[\b a_i;\b b_i]$ is contained in
$\{u=\pm1\}$, and
$(\b b_i,\b a_i)
\setminus\{v=v(\b b_i),v(\b a_i)\}$ is disjoint from $D_i$.
Here
$D_i=\{v(\b a_i)-r_i <v <v(\b b_i)+r_i\}\subset
\{-\frac 12 <v < \frac 12\}$
with
$2r_i=|v(\b b_i)-v(\b a_i)|$ and, finally,

\item[iv)]
the sets $D_i$, $i\in\{1,2,...,n\}$, are pairwise disjoint.

See Figure \ref{fig:72}
 keeping in mind  that $D_i$
is the sub rectangle of $D$ with vertical edges $\{\tilde s=x_i\}$
and $\{\tilde s=x_i'\}$. See also Figure \ref{fig:72} keeping in
mind  that $D_i$ is the rectangle with vertices $\b z_i'$, $\b
z_i$, $\b w_i'$, and $\b w_i$.
\end{itemize}

By lemma \ref{lem:72}, the amount of lifting gained in each set $D_i$ is $2c$
times $r_i$ and it is carried to $[\b a_{ i+1};\b b_{ i+1}]$,
rescaled almost
linearly, by the $\mathcal M$- mean   curvature foliation.
Consequently, as $nc$ is near
1, all of these lifting can be added up as required.

 So,  all recurrent minimal $\mathcal M$- mean   curvature lines can be
eliminated.

To eliminate the recurrent maximal $\mathcal M$- mean curvature
lines of  $\mathbb H_{\alpha,2}^{\mathcal M}$, it is  necessary to perform the same
deformation analysis as  above, applied to this
case, with no fundamental change.

The $\mathcal M$- mean   curvature
separatrices of $\mathbb H_{\alpha,1}^{\mathcal M}$ are stabilized taking care to consider
the
$C^2-$deformations of the immersion $\alpha$,  with support in $\mathcal M$- mean
curvature charts disjoint from  the  nowhere dense set A (see proposition
\ref{prop:72}) consisting  of  the minimal $\mathcal M$- mean
curvature  separatrices and the minimal $\mathcal M$- mean   curvature cycles.
 So the stabilized  minimal $\mathcal M$- mean curvature  separatrices and minimal
$\mathcal M$- mean   curvature cycles
are preserved and these deformations
do not produce any new non trivial
recurrence for the minimal $\mathcal M$- mean   curvature lines.

Therefore the immersion $\alpha $ can be approximated in the $C^2-topology$ by an immersion $\alpha_1$ having all minimal and maximal $\mathcal M$- mean   curvature  separatrices stabilized.

\begin{figure}[htbp]
 \begin{center}
 \includegraphics[angle=0, width=11cm]{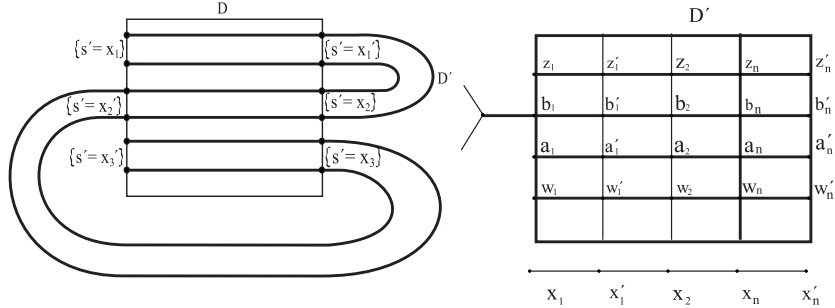}
 \caption{ Lifting of $\mathcal M$- mean   curvature lines  \label{fig:72} }
    \end{center}
  \end{figure}

 \noindent {\bf  Part 2: Conclusion of the Proof of  Theorem \ref{th:d}}

The first step is to approximate an  immersion $\alpha$ of compact and oriented surface
 $\mathbb V^2$ by an immersion having all umbilic points of the type $M_i$, $i=1,2,3$, proposition \ref{prop:2}, and all parabolic points of   types as described in propositions
\ref{prop:p1}, \ref{prop:p2},  \ref{prop:p1b} and \ref{prop:p2b}.

 This can be done by the Transversality Theorem establishing the condition $T=b(b-a)\ne 0$ and by a
 finite number of small local changes on the coefficients of the third jet of the immersion at the
umbilic points. Similar for parabolic points, see \cite{bgm}.

Next approximate the immersion $\alpha$ in $C^s$-sense by an analytic immersion,
 which will be denoted by $\alpha_1$.
There are two cases to consider.

 \noindent {\it Case 1.} The region $\mathbb E \mathbb V^2_{\alpha_1}$ is diffeomorphic to an annulus and the immersion $\alpha_1$ is without umbilics and tangential parabolic singularities.

By using proposition \ref{prop:72} it is possible to obtain an analytic immersion $\alpha^\prime$, $C^2$ close to
$\alpha_1$ having only finitely many $\mathcal M$- mean   curvature cycles, all of which have finite multiplicity.

The resulting immersion $\alpha^\prime$ can be deformed around a   $\mathcal M$- mean curvature  cycle to a obtain an
immersion with a hyperbolic  $\mathcal M$- mean   curvature cycle. If this immersion is approximated  by an analytic
one, $\alpha^{\prime\prime}$, will have only finitely many   $\mathcal M$- mean   curvature cycles, all  of which with
finite multiplicity. In either case, using proposition \ref{prop:42}, $\alpha^{\prime\prime}$ can  be approximated
by an immersion $\tilde{\alpha}$, all whose   $\mathcal M$- mean   curvature cycles are hyperbolic, which  belong to the
class  ${\mathcal G}_i({\mathbb V^2})$, $i=1, 2$,
since conditions i), iii),   iv) and v) are guaranteed by the Stabilization Lemma \ref{lem:74}.
 This ends the proof in this case.

\noindent {\it Case 2.}  The analytic immersion $\alpha_1$ has
separatrices which can be associated to umbilic points, all are of
the types $M_i$, $i=1,2,3$, or parabolic tangential singularities
as in section \ref{sec:5}.

 In this case, as shown in Part 1,  the immersion $\alpha_1$ can be taken so that both, minimal and maximal  $\mathcal M$- mean   curvature   separatrices are stabilized and without non trivial recurrences.
The next step, using proposition \ref{prop:42}, is to deform the immersion in order to obtain an immersion with all minimal and maximal $\mathcal M$- mean   curvature cycles hyperbolic.
  This ends the proof.

\begin{remark}\label{rm:gh}
With the proof of Theorem \ref{th:d} for a fairly general $\mathcal M$, we have completed the density results formulated in \cite{g} and \cite {h} for the particular cases $\mathcal M =\sqrt{ \mathcal K}$ and
 $\mathcal M =\frac{ \mathcal K}{ \mathcal H}$.
\end{remark}

\section{ Concluding Remarks}\label{sec:8}

This paper  presents a theory of unification and generalization for
 the   {\it classical mean curvature} lines on a surface immersed
  in $\mathbb R^3$. See section \ref{sec:1} and the papers
  \cite{a1,   m, g, h, a2}. To this end  the  notion of
   {\it  mean curvature function} was introduced,
    assimilating and adapting for  our   purposes
    in  Geometry the general properties for the {\it Means} studied
    in  Arithmetic and Analysis. See Chapter 8 of Borwein and Borwein \cite{bb}.

The Structural Stability, with its well established achievements
 in the Differential Equations of Geometry,   has  been set as
 a primary goal and a test for  the penetration of the
 generalization proposed in this paper. The methods developed
 in previous specific  works dealing with the {\it classical  means}  have been further
 elaborated and adapted to apply to  the general case of differential equations
 of $\mathcal M$-{\it mean curvature lines} treated here.

It has been established that  the {\it umbilic points} and the
  {\it cycles} of the $\mathcal M$-{mean curvature foliations}
  present  a remarkable  analogy with those of  the classical
   {\it Arithmetic, Geometric} and {\it Harmonic - mean curvature}
   corresponding  cases. See sections \ref{sec:3}  and \ref{sec:4}.

The study of the parabolic singularities  revealed new interesting
aspects. The  analysis presented here was possible by imposing
special regularity conditions  on the  general {\it mean curvature
function} $\, \mathcal M$. Nevertheless, significant cases were
covered. The results achieved, however,  do not apply to the
classic  {\it AG mean} (example \ref {ex:2}), which has a
parabolic pattern not  reducible to algebraic form called {\it
1/k-regularity}. See section \ref{sec:5} for the study of the
cases $k=1,\, 2$.

The development of the
transcendental analysis needed to study the parabolic
singularities of   {\it AG mean} can be regarded as the first
problem of interest left open in this paper.

A great deal  of problems -- such as the studies of   bifurcations and  of immersions of higher dimension and co-dimension --  already proposed in the particular cases of {\it classical mean curvatures},  make sense and have a renewed challenge in the present generalized setting.

By taking  the function $\mathcal M$  as a functional parameter, or itself depending or  a real parmeter, as in the case of the {\it Holder Mean of order r} denoted  ${\mathcal H}_r$ in example \ref {ex:1},
the bifurcation analysis of the transition between different classical differential equations of geometry  and pertinent foliations with singularities   gains new vitality.

The most intriguing of these problems is  the {\it Closing Lemma}. In fact, to prove how to raise the class proximity class from $C^2$  to $C^3$ in Theorem \ref{th:d} is not known even for the case principal foliations. See section \ref{sec:7} and \cite{gsln, gs2}.
Also to achieve the   $C^1$ density for Structural Stability of folded,  \lq\lq billiard\rq\rq,   non-convex configurations is also an open problem in all cases of nets, including {\it asymptotic} ones. See section \ref{sec:6} and \cite{a2}.

\vskip .5cm

\author{\noindent Jorge Sotomayor\\Instituto de Matem\'{a}tica e Estat\'{\i}stica,\\Universidade de S\~{a}o Paulo,
\\Rua do Mat\~{a}o 1010, Cidade Universit\'{a}ria, \\CEP 05508-090, S\~{a}o Paulo, S.P., Brazil \\
\\ Ronaldo Garcia\\Instituto de Matem\'{a}tica e Estat\'{\i}stica,\\
Universidade Federal de Goi\'as,\\CEP 74001-970,
Caixa Postal 131,\\Goi\^ania, GO, Brazil}

\end{document}